\documentclass[11pt,twoside]{article}
\usepackage{mathbbold}

\usepackage{amsmath, amssymb, mathrsfs, graphicx}

\def\az{\alpha}  \def\bz{\beta}
    \def\dz{\delta}
\def\ez{\eta}    \def\fz{\varphi}
\def\gz{\gamma}  
\def\lz{\lambda} \def\mz{\mu}

\def\rz{\rho}        
        \def\uz{\theta}
\def\vz{\varepsilon}

\def\qd{\quad}
\def\qqd{\qquad}

\newcommand{\mathsym}[1]{{}}

\newcommand{\cirnm}[1]{\bigcirc{\mbox{\hspace{-0.75em}}#1}\,}

\def\le{\leqslant}
\def\ge{\geqslant}

\font\cms=cmss9 scaled \magstep1

\def\nnd{\noindent}

\def\crl{\nnd\bg{crl1}}
\def\lmm{\nnd\bg{lmm1}}
\def\prp{\nnd\bg{prp1}}

\def\xmp{\nnd\bg{xmp1}}

\def\decrl{\end{crl1}}
\def\delmm{\end{lmm1}}
\def\deprp{\end{prp1}}
\def\dexmp{\end{xmp1}}

\def\prf{\medskip \noindent {\bf Proof}. }
\def\qed{\text{\quad $\Box$}}
\def\deprf{\qed\medskip}
\def\bg{\begin}
\def\be{\bg{equation}}
\def\de{\end{equation}}

\def\dear{\end{eqnarray}}
\def\lb{\label}
\def\ct{\cite}

\newcommand{\rf}[2]{[\ref{#1}; #2]}

\def\den{\end{enumerate}}

\def\bbb{\mathbbold}

\begin{document}


\thispagestyle{empty}
\renewcommand{\thefootnote}{\fnsymbol{footnote}}

\vspace*{-0.9in}
\noindent {Front. Math. China 2016, 11(6): 1379--1418}\newline
\noindent {DOI 10.1007/s11464-016-0573-4}

\begin{center}
{\bf\Large Efficient initials for computing maximal eigenpair}
\vskip.15in {Mu-Fa Chen}
\end{center}
\begin{center} (Beijing Normal University)\\
\vskip.1in June 21, 2015
\end{center}
\vskip.1in

\markboth{\sc Mu-Fa Chen}{\sc Computing the maximal eigenpair}



\date{}


\footnotetext{Received October 22, 2015; accepted June 23, 2016.}

\begin{abstract} This paper introduces some efficient initials for a well-known
algorithm (an inverse iteration) for computing the maximal eigenpair of a class of real matrices. The initials not only avoid the collapse of the algorithm but are also unexpectedly efficient.
The initials presented here are based on our analytic estimates of the
maximal eigenvalue and a mimic of its eigenvector
for many years of accumulation in the study of stochastic stability speed.
In parallel, the same problem for computing the next to the maximal eigenpair is also studied.
\end{abstract}

\nnd {\small 2000 {\it Mathematics Subject Classification}: 15A18, 65F15, 93E15, 60J27}
\medskip

\nnd {\small {\it Key words and phrases}. Perron-Frobenius theorem, power iteration, Rayleigh quotient iteration, efficient initial, tridiagonal matrix, $Q$-matrix.}

\bigskip

\section{Introduction. Two algorithms and a typical\\ example}

   Consider a nonnegative irreducible matrix $A=(a_{ij})$ on $E:=\{0, 1, \ldots, N\}$, $N<\infty$. By the well-known
Perron-Frobenius theorem, the matrix has uniquely a positive eigenvalue $\rho(A)$ having positive left-eigenvector and positive right-eigenvector. Moreover, both the left-eigenspace and the right-eigenspace of $\rho(A)$ have dimension one. This eigenvalue is maximal
in the sense that for every other eigenvalue $\lz_k$, we have $\rho(A)\ge |\lz_k|$.
The last equality sign appears only if $A$ has a period $p>1$. For instance, for
$$A=\begin{pmatrix}
0 & 0 & 1\\
0 & 0 & 1\\
1 & 1 & 0\\
\end{pmatrix},$$
we have $p=2$ and the eigenvalues of $A$ are $\pm \sqrt{2}$ and $0$.
However, we may assume that $\rho(A)> |\lz_k|$ for every other eigenvalue $\lz_k$.
Actually, if $\lz=\rho e^{i\uz}$ with $\uz\ne k\pi/2$ for every odd $k\in {\mathbb Z}$, then
for every $\vz>0$, we have $\rho+\vz> |\rho e^{i\uz}+\vz|$. This means the required assertion holds for the shifted pair $\rho+\vz$ and $\lz+\vz$. In other words,
an analog of the Perron-Frobenius
theorem is meaningful for the matrices having nonnegative off-diagonal elements only,
their diagonal elements can be arbitrary but real. By a shift if necessary, such
a matrix can be transformed into a nonnegative one: the maximal eigenvector is kept
but their maximal eigenvalues are shifted from one to the other.
In this paper, we are interested in computing $\rho(A)$ and its corresponding eigenvector. This is a very important problem, we will come back to its motivation in the next section.

There are mainly two popular algorithms for this problem. Unless otherwise stated,
the eigenvector below means the right-eigenvector. Then, the maximal eigenpair (the maximal eigenvalue and its eigenvector) is denoted by $(\rho(A), g)$.

\medskip
\nnd{\bf Power iteration}.\quad{\cms Given an initial vector $v_0\in {\mathbb R}^{N+1}$
having a nonzero component in the direction of $g$ with $\|v_0\|=1$, define
\be v_{k}=\frac{A v_{k-1}}{\|A v_{k-1}\|},\quad z_k={\|Av_k\|},\qqd k\ge 1,\de
where $\|\cdot\|$ is an arbitrary but fixed vector norm.
Then $v_k$ converges to the eigenvector $g$ of $\rho(A)$ and $z_k$ converges to $\rho(A)$ as $k\to\infty$.}\medskip

Even it is not necessary, in the next algorithm, we fix the vector norm to be the Euclidean one (or equivalently, the $\ell^2$-norm). Actually, a
refined choice is using the inner product and the norm in the space $L^2(\mu)$ for a suitable measure $\mu$ to be specified case by case, as illustrated by the improved algorithm given at the end of
Sections \ref{s-03} and \ref{s-04}. See also Section \ref{s-06}.

\medskip
\nnd{\bf Rayleigh quotient iteration (a variation of inverse iteration)}.\quad{\cms Choose a pair $(z_0, v_0)$
as  an approximation of $(\rho(A), g)$ with $v_0^*v_0=1$, where $v^*$ is the transpose of $v$. In particular, one may set $z_0=v_0^*Av_0$ for a given $v_0$. At the $k$th $(k\ge 1)$ step, solve the linear equation in $w_k$:
\be (A-z_{k-1}I) w_k =v_{k-1},  \lb{02}\de
where $I$ is the identity matrix on $E$, and define
$$v_k=\frac{w_{k}}{\sqrt{w_{k}^*w_k}},\qqd z_k=v_{k}^*Av_k.$$
If the pair $(z_0, v_0)$ is close  enough to $(\rho(A), g)$, then $(z_k, v_k)$ converges to $(\rho(A), g)$ as $k\to\infty$.
}\medskip

In what follows, unless otherwise stated, we fix $z_0$ to be the particular
choice just defined. We now use a typical example (which will be studied time by time in the paper) to illustrate the effectiveness and their
difference of the above two algorithms.

\xmp\lb{t-01} \;\;{\cms Let $E=\{0, 1, \ldots, 7\}$ and
$$Q=\left(
\begin{array}{cccccccc}
 -1 & 1 & 0 & 0 & 0 & 0 & 0 & 0 \\
 1 & -5 & 2^2 & 0 & 0 & 0 & 0 & 0 \\
 0 & 2^2 & -13 & 3^2 & 0 & 0 & 0 & 0 \\
 0 & 0 & 3^2 & -25 & 4^2 & 0 & 0 & 0 \\
 0 & 0 & 0 & 4^2 & -41 & 5^2 & 0 & 0 \\
 0 & 0 & 0 & 0 & 5^2 & -61 & 6^2 & 0 \\
 0 & 0 & 0 & 0 & 0 & 6^2 & -85 & 7^2 \\
 0 & 0 & 0 & 0 & 0 & 0 & 7^2 & -113
\end{array}
\right).$$
Then we have $\rho(Q)\approx -0.525268$ with eigenvector
$$\approx({55.878,\; 26.5271,\; 15.7059,\; 9.97983,\; 6.43129,\; 4.0251,\; 2.2954,\; 1})^*.$$
Starting from $v_0$ which is the normalized vector of
$$(1,\; 0.587624,\; 0.426178,\; 0.329975,\; 0.260701,\; 0.204394,\; 0.153593, \;
0.101142)^*,$$
the power iteration (applied to the nonnegative $A:=113\, I +Q$) arrives at $-0.525268\approx \rho(Q)$
after 990 iterations. Here we adopt the $\ell^1$-norm:
$$\|v\|=\sum_{k\in E} |v_k|.$$
}\dexmp

We now give a little more details about the computations for this example.
\vspace{0.05truecm}
\def\tba{\begin{array}{cc}
 0 & 2.11289 \\
 1 & 1.42407 \\
 2 & 1.37537 \\
 3 & 1.22712 \\
 4 & 1.1711 \\
 5 & 1.10933 \\
 6 & 1.06711 \\
 7 & 1.02949 \\
 8 & 0.998685 \\
 9 & 0.971749 \\
 10 & 0.948331 \\
 11 & 0.927544 \\
 12 & 0.908975 \\
 13 & 0.892223
 \end{array}}

 \def\tbb{\begin{array}{cc}
 14 & 0.877012 \\
 15 & 0.86311 \\
 16 & 0.850338 \\
 17 & 0.838548 \\
 18 & 0.827619 \\
 19 & 0.817449 \\
 20 & 0.807953 \\
 30 & 0.738257 \\
 40 & 0.694746 \\
 50 & 0.664453 \\
 60 & 0.641946 \\
 70 & 0.624473 \\
  80 & 0.610468 \\
 90 & 0.598963
 \end{array}}

 \def\tbc{\begin{array}{cc}
 100 & 0.589332 \\
 120 & 0.574136 \\
 140 & 0.56279 \\
 160 & 0.554157 \\
 180 & 0.547529 \\
 200 & 0.542423 \\
 300 & 0.529909 \\
 400 & 0.526517 \\
 500 & 0.525603 \\
 600 & 0.525358 \\
 700 & 0.525292 \\
 800 & 0.525274 \\
 900 & 0.52527 \\
 \ge 990 & 0.525268
\end{array}}
\vspace{-0.8truecm}

\begin{center}Table 1.\quad Partial outputs of $(k, -z_k)$.\end{center}
\vspace{-0.6truecm}
\begin{center}
  \begin{tabular}{| c | c | c | }
    \hline
$\tba$  & $\tbb $ & $\tbc$\\ \hline
  \end{tabular}\vspace{-0.25truecm}
 \end{center}

\nnd Table 1 gives us partial outputs of $(k, -z_k)$. The corresponding figure below shows that $-z_k$ decreases quickly for small $k$, but the convergence goes very slow for large $k$.
\vspace{-0.4truecm}

{\begin{center}{\includegraphics[width=12.5cm,height=7.0cm]{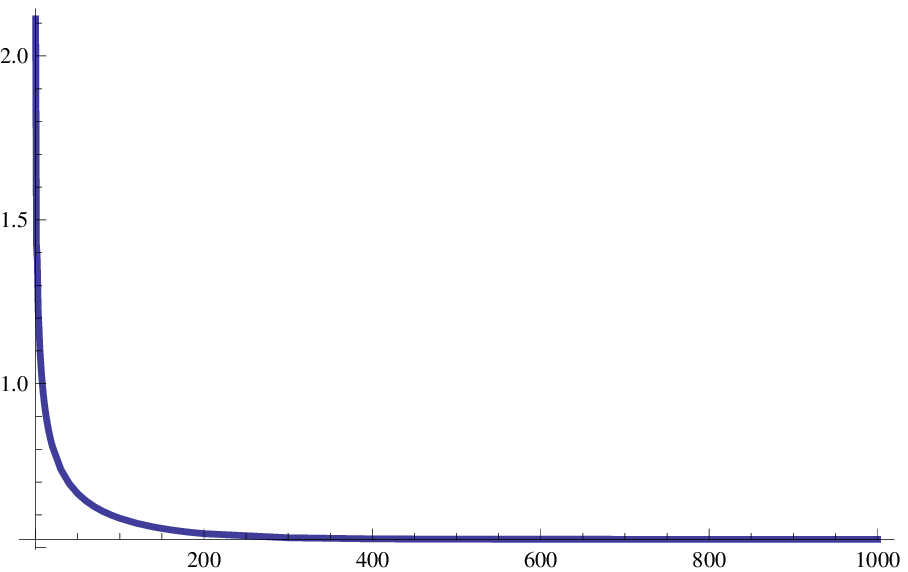}
\vspace{-0.3truecm}}\end{center}}

\begin{center}{\bf Figure 1}\lb{fig-1}\quad The figure of $-z_k$ for $k=0, 1, \ldots, 1000$.\end{center}

The advantage of the first algorithm is that it allows us to use a quite
arbitrary positive initial vector $v_0$. The reason why the convergence of the example at the beginning steps goes quite fast is because we have used a very
good initial $v_0$, as will be studied in Section \ref{s-03}. However, for larger $k$, the convergence becomes very
slow, that is the limitation of this algorithm. From Figure \ref{fig-1}, it is
clear that one may stop the computation at 300 iterations since then
the results are almost the same. However, we keep going on until the
six precisely significant digits as limited by a computer using
Mathematica 9. The reason for doing so is for the comparison with
other algorithms to be studied later.

Certainly, we expect the second algorithm to be more efficient.
Now, what can we expect? Since this problem is often used in practice,
we would be very happy if a new algorithm can reduce the number of
iterations seriously, say, 500 for instance. Certainly, we would be
very surprising if it can be reduced to 250. Let us think this question more carefully.
Suppose that we are now interested in the maximal eigenvalue only, and
suppose we know it is located on $(0, 1)$ (actually, as we will see by Proposition \ref{t-11} below, the maximal eigenvalue of $A:=113\,I + Q$ is located in $(0, 113)$). We may use the Golden Section Search (a famous method in optimization theory), its speed is about
$0.618^{-1}$. Then, to obtain the six precisely significant digits as in the
last example, one needs at least 24 iterations since $10^{-6}\approx 0.618^{24}$. Suppose we can adopt a faster algorithm, the Bisection Method, then it requires about 20 iterations since $10^{-6}\approx 2^{-20}$.
Hence, it is reasonable if an algorithm uses more than 20 iterations to arrive at the same precise level. Having this analysis in mind, we were shocked when the next result came to us.

\xmp \lb{t-02} \;\;{\cms The matrix $Q$ and the initial vector $v_0$ are the
same as in the last example but we now adopt the $\ell^2$-norm. The Rayleigh quotient iteration
(applied to $Q$) starts at $z_0=v_0^*Qv_0\approx -0.78458$ and then arrives at the same result as in the last example at the second step
$$z_1\approx -0.528215,\;\; z_2\approx -0.525268.$$
}\dexmp

Example \ref{t-02} is the main illustrating example (which will be further improved by Example \ref{t-07} below) of this paper. It shows
that the second algorithm can be extremely powerful. The key to this result is that we have chosen an efficient initial vector $v_0$ and then the resulting $z_0$ is also close to $\rho(Q)$. It may be the position
to compare the use of $\ell^1$-norm and $\ell^2$-one. Let everything be the
same as in the last example but replacing the $\ell^2$-norm by the $\ell^1$-one. Then the iteration starts at $z_0\approx -0.367937$ and arrives at the same result at the third step:
$$z_1\approx -0.509272,\qd z_2\approx -0.52509,\qd z_3\approx -0.525268.$$
The result comes with no surprising: it is easier to use the $\ell^1$-norm
in the computation but it is a little less efficient than using the  $\ell^2$-norm.

We have seen the power of the second algorithm. However,
``too good'' is dangerous.
Each eigenvalue $\lz_k\ne \rho(A)$ can be a pitfall of the algorithm provided either $z_0$ is close  enough to $\lz_k$ or $v_0$ is close  enough to the eigenvector $g_k$ of $\lz_k$.  The next example
illustrates the latter situation. For which a simpler $v_0$ deduces
its corresponding $z_0$ to be more close to $\lz_2$ rather than
$\lz_3$.

Here and in what follows, we often use the so-called $Q$-matrix
$Q=(q_{ij}: i, j\in E)$ which means that $q_{ij}\ge 0$ for every pair
$i\ne j$ and $\sum_{j\in E}q_{ij}\le 0$ for every $i\in E$. This implies
the intrinsic use of probabilistic idea. For convenience, we often write
by
$$0<\lz_0< |\lz_1| \le |\lz_2| \le\ldots,$$
where $\{\lz_j\}$ is the sequence
of the eigenvalues of $-Q$. Then $\lz_0=-\rho(Q)$.

\xmp \lb{t-03} \;\;{\cms The matrix $Q$ is the same as the last example and we use again the $\ell^2$-norm. Replace $Q$ by $-Q$ (then the corresponding $z_k>0$). Choose the initial vector $v_0$ to be the normalized uniform vector:
$$v_0=\{1, 1, 1, 1, 1, 1, 1, 1\}/\sqrt{8}.$$
Then with the particular choice given in the algorithm
$$z_0=v_0^*(-Q)v_0=8,$$
we obtain the following output at the first 4 steps of the iterations:
$$z_1\approx 4.78557,\;\; z_2\approx 5.67061,\;\;
z_3\approx 5.91766,\;\; z_4\approx 5.91867\approx \lz_2.$$
The first two eigenvalues of $-Q$ are
$$\lz_0\approx 0.525268,\qd
\lz_1\approx 2.00758,\qd\lz_3\approx 13.709,$$ respectively. Hence the limit $\lz_2$ is quite away from what we are interested in.}\dexmp

By the way, let us mention that in practice,
we can stop our computation once the components of
the first output $v_1$ have different signs, and try to choose a new
initial pair $(v_0, z_0)$. This is due to the fact that the maximal
vector should be positive/negative up to a constant. Here, in the last example,
$$\aligned
v_1=(&-0.26762,\; 0.242432,\; -0.522646,\; -0.579319,\; -0.423469,\; -0.253452, \\
&-0.124365,\; -0.0425044)^*.\endaligned$$
Each of the components is negative except the second one.

The next example shows that we can still arrive at the expected result
for a good initial $z_0$ even if $v_0$ is quite rough.

\xmp \lb{t-04} \;\;{\cms Everything is the same as in the last example
except
$$z_0=2.05768^{-1}\approx 0.485985.$$
Then $\{z_k\}$ approaches to $\lz_0$ at the second step.
$$z_1\approx 0.525998,\;\;
z_2\approx 0.525268.$$
}\dexmp

The paper is organized as follows. In the next section, we first review the
five sources of the motivation for our problem. Then we recall the known
convergence of these algorithms. From the above examples, we have seen that
the second algorithm is much more attractive. To which, we need a careful
design in choosing the initial pair $(v_0, z_0)$. Clearly, an efficient
initial pair is just a good estimate of the pair in advance.
This itself is a hard topic in the study of eigenvalue problem and so it is
understandable that the initial problem is still largely open
in the eigenvalue computation theory.
A complete, analytic (explicit) solution to this problem is presented in Section \ref{s-03} first for tridiagonal
matrices (after a suitable relabeling if necessary), and then for a class of more general matrices in terms of the
so-called Lanczos tridiagonalization procedure. The main extension to the
general situation is presented in Section \ref{s-04} which consists of two
subsections. In the first one, we concentrated on the construction of
$z_0$ for a fixed simplest $v_0$. The second one in even more technical,
in which we are mainly working on the construction of $v_0$. A number of
examples are illustrated, case by case, for the results in the paper.
It is remarkable that only the one-step iteration scheme, as illustrated by
the two algorithms used above, is used in the paper. In Section \ref{s-05},
we make either additional proofs of some results in the main context,
or additional remarks on related problems.
In particular, we prove a convergence result of our approximating procedure
for the principal eigenvalue of birth--death processes which have been
studied for a long period up to now. A summary for the use of the algorithms up
to Section \ref{s-04} is given at the end of Section \ref{s-05}. The study on the
next eigenpair is delayed to Section \ref{s-06}.

\section{Motivation of problem and convergence of \\
algorithms}

In this section, for the reader's convenience, we recall briefly the motivation
of our problem and the well-known convergence of the two algorithms introduced
in the first section.

\subsection{Motivation}

It seems not necessary to mention the value of the study on the problem
since the matrix eigenvalue computation is used almost everywhere. The next five sources reflect more or less
the road where we started and finally arrive here.

\subsection*{Google's PageRank}

When we search an expression from the network, a large number of
webpages are collected. The question is how to output them on the
screen of our computer. For this, we need to rank the pages. The procedure
goes as follows. According to the connections of the websites, we get a
nonnegative matrix $A$. To which we have the largest eigenvalue $\rho(A)$
and its corresponding po\-sitive left-eigenvalue. The normalized left-eigenvector
gives us an order of the webpages, that is the PageRank as we required.
Nowadays, there are a large number of publications on Google's PageRank,
see for instance \ct{lan1}, in which the power iteration is studied but
not the inverse iteration.

\subsection*{Global optimization of planned economy}

Regarding the matrix $A$ as a structure matrix in economy, Hua \ct{hua84}
proved that the optimal input of the planned economy is the left-eigenvector
$u$ of $\rho(A)$. Surprisingly, Hua \ct{hua84} also proved that if one uses a different
input rather than $u$, then the economy will go to collapse (i.e., some
components of the product in the economic system will become less or
equal to zero). Mathematically, this situation is very much the same as
the last one, but in a completely different context. As far as I know, the
practical algorithms for $(u, \rz(A))$ were not studied carefully during
that period, except a formula was mentioned in \ct{hua84}:
$$\rho(A)=\lim_{\ell \to\infty}\bigg(\frac{\text{\rm Trace ($A^{\ell}$)}}{N+1}\bigg)^{1/{\ell}}.$$

\subsection*{Stationary distribution of time-discrete Markov chain}

If $A$ itself is a transition probability matrix, then the left-eigenvector corresponding
to the largest eigenvalue one is nothing but the stationary distribution
of the corresponding Markov chain. This explains the stability meaning in the two
situations just discussed above. Based on this idea, we obtained a probabilistic
proof of Hua's collapse theorem. Refer also to \rf{cmf05}{Chapter 10}
for additional story and related references.

Computing the stationary distribution of a given Markov chain is
very important in practice and so has been studied quite a lot in the
past years, including the so-called Markov Chain Monte Carlo (MCMC),
perfect/backward coupling, and so on.

\subsection*{Exponential decay of time-continuous Markov chain}

The maximal eigenvalue $\rho(Q)$ in Example \ref{t-01} describes the
exponential decay rate of the Markov chain with semigroup $(P_t= e^{tQ}: t\ge 0)$. The present paper is based on our study on this topic, as will be
seen from the subsequent sections.

\subsection*{Phase transitions}

The last topic and the investigation on related stability speed are actually motivated from the study on phase transitions in statistical mechanics
(cf. \ct{cmf04, cmf05} for more references within). This is a challenge
topic in mathematics since it is mainly in infinite-dimensional setting.
To which, the mathematical tools are rather limited. Therefore we have to look for new tools or develop some known traditional tools. To this end, we have already
visited several branches of mathematics, including the computation theory. We are now glad to be able to say something on the last field
after a long trip of the study.
\medskip

In the second part of this section, we review some well-known facts
on the convergence of the algorithms.

\subsection{Convergence of the algorithms}

Here is the convergence of the power iteration. In this subsection, we
suppose that the given matrix $A$ (not necessarily nonnegative) has the
dominant eigenvalue $\lz_0$
(i.e.,  $|\lz_0|> |\lz_j|$ for all other eigenvalues $\lz_j$)
which is simple. The extension to the periodic situation is also
possible, but is omitted here, one simply replaces the convergence of the original
sequence by a subsequence.

\lmm\lb{t-05}\;\;{\cms Suppose that the initial vector $v_0$ has a nonzero component in the direction of the dominant eigenvector $g$. Then
$$v_{k}=\frac{A^k v_0}{\|A^k v_0\|}\to g\;\;\text{\cms and}\;\;
v_k^* A v_k \to \lz_0\qqd \text{\cms as }k\to\infty.
$$
Moreover,
$$\lim_{n\to\infty}\frac{A^n v_0}{A^{n-1} v_0}=\lz_0,$$
where for given vectors $u$ and $v$, the ratio $u/v$ is understood as the
quotient function of the functions $u$ and $v$.}
\delmm

\prf Suppose that the eigenvalues are all different for simplicity. Otherwise, one simply uses the Jordan representation of matrices. Write
$$v_0=\sum_{j=0}^N c_j g_j$$
for some constants $(c_j)$ with $g_0=g$. Then $c_0\ne 0$ by assumption and
$$A^k v_0=\sum_{j=0}^N c_j\lz_j^k g_j
=c_0\lz_0^k \bigg[g+\sum_{j=1}^N \frac{c_j}{c_0}\bigg(\frac{\lz_j}{\lz_0}\bigg)^k g_j\bigg].
$$
Since $|\lz_j/\lz_0|<1$ for each $j\ge 1$ and $\|g\|=1$, we have
$$\frac{A^k v_0}{\|A^k v_0\|}\to \frac{c_0}{|c_0|} g\qqd \text{as }k\to \infty,$$
and then
$$v_k^* A v_k\to g^* A g=g^* \lz_0 g =\lz_0 \qqd \text{as }k\to \infty.$$
We have thus proved the main assertion of the lemma. The proof of the last
assertion is similar. \deprf

\smallskip

Clearly, the convergence speed in the lemma is
$$\bigg|\frac{\lz_1}{\lz_0}\bigg|^k, \qqd |\lz_1|:=\max\{|\lz_j|: j>0\}.$$

The next result is the convergence for the inverse iteration.

\lmm\lb{t-06}\;\;{\cms Under the assumption of the last lemma, for each $z$ close  to $\lz_0$, we have
$$v_{k}=\frac{(A-zI)^{-k} v_0}{\|(A-zI)^{-k} v_0\|}\to g\;\;\text{\cms and}\;\;
v_k^* A v_k \to \lz_0\qqd \text{\cms as }k\to\infty.
$$
Moreover,
$$\lim_{n\to\infty}\frac{(A-zI)^{-n} v_0}{(A-zI)^{-n+1} v_0}=\frac{1}{\lz_0-z}.$$
}\delmm

\prf Note that for $z$ close to $\lz_0$, the dominant eigenvalue of the
matrix $(A-z I)^{-1}$ is $(\lz_0-z)^{-1}$ with the same dominant eigenvector
$g$ as the one for $A$. The proof is very much the same as the previous one.
$\lz_0$
\deprf

The iteration given in Lemma \ref{t-06} is called the inverse iteration. It is remarkable that the convergence speed in this lemma is
$$\bigg|\frac{\lz_0-z}{\lz_1-z}\bigg|^k\sim \bigg|\frac{\lz_0-z}{\lz_1-\lz_0}\bigg|^k$$
when $z$ is sufficiently close to $\lz_0$.
At the $k$th step, replacing $z$ by the Rayleigh quotient approximation
$z_k=v_{k}^*Av_k$, we obtain the Rayleigh quotient iteration as described
in the first section.  Clearly, the last algorithm is an acceleration of
the inverse iteration.
The price is that the initial $z_0$ has to be chosen close to $\lz_0$ which is usually not
explicitly known. Otherwise, if $z_0$ is chosen close to some $\lz_j\ne \lz_0$, then a similar proof of Lemma \ref{t-06} shows that $v_k^* A v_k$ converges to the pitfall $\lz_j\ne\lz_0$.
In practice, once $z=z_0$ is chosen in a suitable neighborhood of $\lz_0$, the
sequence $z=z_k$ converges to $\lz_0$ rapidly, as illustrated by Examples
\ref{t-02} and \ref{t-04}.
More precisely, Example \ref{t-01} applies the power iteration to $A:=113 I + Q$,
its convergence speed is
$$\sim\bigg(\frac{113-\lz_1}{113-\lz_0}\bigg)^k\approx \bigg(\frac{113-2.00758}{113-0.525268}\bigg)^k\qqd \text{as } k\to \infty.$$
Examples \ref{t-02} and \ref{t-04} use the Rayleight quotient iteration which has the
convergence speed
$$\sim \prod_{j=0}^k \frac{\lz_0-z_j}{\lz_1-z_j}\qqd \text{as } k\to \infty.$$
Since $z_k\to \lz_0$, the last convergence is much fast than the previous one.
Honestly, this still does not
answer the reason why the inverse algorithm in Example \ref{t-02} can achieve the six
significant digits at the second iteration.

\section{Efficient initials. Tridiagonal case}\lb{s-03}

Again, assume that $A=(a_{ij})$ on $E=\{0,1,\ldots N\}$, $N<\infty$ is
irreducible and having non-negative off-diagonal elements. Assume also
that the matrix is tridiagonal (after a suitable relabeling if necessary): $a_{ij}=0$ unless $|i-j|\le 1$. By a shift
$Q:=A-m I$
 if necessary, where $I$ is the identity matrix on $E$ and
 $m=\max_{i\in E}\sum_{j\in E}a_{ij},$
 one may assume that
$$Q\!=\!\left(\!\begin{array}{ccccc}
-(b_0+c_0)\!\!\! & b_0 &0&0 &\cdots \\
a_1 &\!\!\! -(a_1 + b_1+c_1) & b_1 &0 &\cdots \\
0& a_2 &\!\!\! -(a_2 + b_2+c_2) & b_2 &\cdots \\
\vdots &\vdots &\ddots &\ddots &\ddots \\
0& 0& \qquad 0 &\quad a_N^{}\;\; & -(a_N^{}+c_N^{})
\end{array}\!\right)\!,$$
where $a_i> 0,\; b_i> 0,\; c_i\ge 0$ but $c_i\not\equiv 0$.
Define
$$\gathered
\mu_0=1, \;\; \mu_n=\mu_{n-1}\frac{b_{n-1}}{a_n}=\frac{b_0 b_1 \cdots b_{n-1}}{a_1 a_2 \cdots a_n},
\qquad 1\le n\le N.
\endgathered
$$
We now split our discussion into two cases.

{\it Case $1$}.\qd Let
$$c_0=\cdots=c_{N-1}=0.$$
Then we may assume that
$c_N>0$. Otherwise, $Q$ has the trivial maximal eigenvalue $0$ with
eigenvector with components being one everywhere. In this case, we
rewrite $c_N$ as $b_N$, ignoring the sequence $(c_i)$, and define
\be\fz_n=\sum_{k= n}^N \frac{1}{\mu_k b_k}, \qqd 0\le n\le N.\lb{3-0}\de

{\it Case $2$}.\qd Let some of $c_i\, (i=0, 1, \ldots, N-1)$ be positive.
Then, we need more work. Define
$$\gathered
r_0=1+\frac{c_0}{b_0},\;\; r_n=1+\frac{a_n+c_n}{b_n}-\frac{a_n}{b_n r_{n-1}},\qqd 1\le n<N,\\
h_0=1,\;\; h_n=h_{n-1}r_{n-1}=\prod_{k=0}^{n-1}r_k,\qqd 1\le n \le N,
\endgathered
$$
and additionally,
$$h_{N+1}= c_N h_N+ a_N (h_{N-1}-h_N).
$$
Finally, define
\be \fz_n=\sum_{k= n}^N \frac{1}{h_k h_{k+1}\mu_k b_k}, \qqd 0\le n\le N
\lb{3-0-1}\de
with a convention that $b_N=1$ to save our notation.

 We remark that in the special case that
 $$c_0=\cdots=c_{N-1}=0,$$
 by induction, it is easy to check that
 $$r_0=\cdots=r_{N-1}=1$$
 and hence
 $$h_0=\cdots=h_{N}=1.$$
 Furthermore, $h_{N+1}=c_N$. Thus, once replacing $c_N$ by $b_N$,
 we return to (\ref{3-0}) from (\ref{3-0-1}).

To state our algorithm, we need one more quantity.
$$\dz_1=\max_{0\le n\le N} \bigg[\sqrt{\fz_n}\sum_{k=0}^n \mu_kh_k^2 \sqrt{\fz_k}+ \frac{1}{\sqrt{\fz_n}}\sum_{n+1\le j \le N}\mu_jh_j^2\fz_j^{3/2}\bigg].  $$

\nnd{\bf Rayleigh quotient iteration in the tridiagonal case}.\quad{\cms For given tridia\-gonal
matrix $A$, define $m$, $(a_i, b_i, c_i)$, $(h_i)$, $(\fz_i)$ and $\dz_1$ as above. Set
$${\tilde v}_0(i)=h_i \sqrt{\fz_i},\qqd 0\le i \le N,\qqd
v_0=\frac{{\tilde v}_0}{\sqrt{{\tilde v}_0^*{\tilde v}_0}}, \qqd z_0=\frac{1}{\dz_1}.$$
At the $k$th step $(k\ge 1)$, solve the linear equation in $w_k$:
\be (-Q - z_{k-1}I) w_k =v_{k-1} \lb{3-1}\de
and define
$$v_k=\frac{w_{k}}{\sqrt{w_{k}^*w_k}},\qqd z_k=v_{k}^*(-Q)v_k.$$
Then $v_k$ converges to $g$ and $m-z_k$ converges to $\rho(A)$ as $k\to\infty$.
}\medskip

It is an essential point that the choice of $z_0$ avoids the collapse since we have known that
$\lz_0(Q)=\lz_{\min}(-Q)\; (\text{the minimal eigenvalue of } -Q)\ge \dz_1^{-1}$ by \rf{cmf10}{Corollary 3.3}.
As an application of this result to Example \ref{t-01}, we have $c_i\equiv 0$
but $b_7=64$ and then $h_i\equiv 1$. We can define $\fz$ by (\ref{3-0}) and then ${\tilde v}_0=\sqrt{\fz}$ which
is the one, up to a free factor $\sqrt{\fz_0}$,
used in Example \ref{t-01}. This is the meaning of ``very good'' claimed in the first section.
We now compute the minimal eigenvalue of $-Q$ using not only ${\tilde v}_0$ but also $\dz_1$.

\xmp \lb{t-07} \;\;{\cms The matrix $Q$ and the vector ${\tilde v}_0$ are the same as in Example $\ref{t-01}$:
$$(1,\; 0.587624,\; 0.426178,\; 0.329975,\; 0.260701,\; 0.204394,\; 0.153593, \;
0.101142)^*.$$
We have $\dz_1=2.05768.$ Then, with the new
$z_0:=\dz_1^{-1}\approx 0.485985$, the Rayleigh quotient iteration
arrives at the expected estimate at the second step:
$$z_1\approx 0.525313,\;\; z_2\approx 0.525268.$$
}\dexmp

Comparing the approximation value of $z_1$ here and that in Example \ref{t-02}, it is clear that
this result is sharper than Example \ref{t-02} (see also the comment below Corollary \ref{t-11-1}).

Now, let us discuss the effectiveness of our algorithm
with respect to the size $N$ of the matrix. In computational mathematics,
one often expects the number of iterations $M$ grows up no more than $N^{\az}$ for
some $\az>0$. It is unusual if $M\approx \log N$ for large $N$. To this
question, considering the basic Example $\ref{t-01}$ with varying $N$, the answer given below is worked out by Yue-Shuang Li using the algorithm introduced in this section and the software MatLab on a notebook. In the first line of Table 2, the reason we use $N+1$ rather than $N$ is that the space is labeled starting at $0$ but not $1$.
\begin{center}Table $2$.\quad  For different $N$, the
eigenvalue $\lz_0$, its lower bound $\dz_1^{-1}$ and $z_1, z_2$.
\end{center}
\vskip-0.4truecm
$$
\hfil\vbox{\hbox{\vbox{\offinterlineskip
 \halign{&\vrule#&\strut\quad\hfil#\hfil\quad&\vrule#&
 \quad\hfil#\hfil\quad&\vrule#&
 \quad\hfil#\hfil\quad\cr
    \noalign{\hrule}
   height2pt&\omit&&\omit&&\omit&&\omit&\cr
  & {$\pmb{N+1}$} &&${\pmb{z_0=\dz_1^{-1}}}$ &&${\pmb{z_1}}$ &&${\pmb{z_2=\lz_0}}$ &\cr
   height2pt&\omit&&\omit&&\omit&&\omit&\cr
   \noalign{\hrule}
   height2pt&\omit&&\omit&&\omit&&\omit&\cr
  & \!100\! &&\!0.348549\!  && \!0.376437\! && {0.376383}  &\cr
   height2pt&\omit&&\omit&&\omit&&\omit&\cr
   \noalign{\hrule}
   height2pt&\omit&&\omit&&\omit&&\omit&\cr
& 500 &&0.310195  && 0.338402 && 0.338329  &\cr
   height2pt&\omit&&\omit&&\omit&&\omit&\cr
   \noalign{\hrule}
   height2pt&\omit&&\omit&&\omit&&\omit&\cr
& 1000 &&0.299089 && 0.32732 && 0.32724  &\cr
   height2pt&\omit&&\omit&&\omit&&\omit&\cr
   \noalign{\hrule}
   height2pt&\omit&&\omit&&\omit&&\omit&\cr
& 5000 &&0.281156 && 0.308623 && 0.308529  &\cr
     height2pt&\omit&&\omit&&\omit&&\omit&&\omit&\cr
   \noalign{\hrule}
      height2pt&\omit&&\omit&&\omit&&\omit&&\omit&\cr
& 7500 &&0.277865 && 0.305016 && 0.304918  &\cr
     height2pt&\omit&&\omit&&\omit&&\omit&&\omit&\cr
   \noalign{\hrule}
& $10^4$ &&0.275762 && 0.30266 && 0.302561  &\cr
 \noalign{\hrule}}}}}\hfill$$
Is it believable? Yes, we have justified the outputs in two different ways:
in each case, firstly, the outputs starting from $z_2$ become the same
(which actually coincides with the output of $\lz_0$). Secondly, by using $v_2$, we can
find upper/lower estimates $\bar\xi$/$\underline\xi$ of $\lz_0$ such that
$z_2\in \big({\underline\xi}, {\bar\xi}\big)$ and moreover
$${\bar\xi}\big/{\underline\xi}\approx 1+ 10^{-5}.$$

The next example is due to L.K. Hua \ct{hua84} in the study of economic optimization (cf. \rf{cmf05}{Chapter 10}). Note that here we are studying the right-eigenvector, the matrix $A$ used below is the transpose of the original one.

\xmp\lb{t-08} \;\;{\cms Let
$$A=\frac{1}{100}\left(\begin{matrix}25 & 40\\
14 & 12\end{matrix}\right).$$
Then
$$\rho (A)=\big(37 + \sqrt{2409}\big)/200\approx 0.430408.$$
With the initials:
$$v_0\approx (0.429166,\; 0.220573)^*\qd\text{\cms and}\qd
z_0:=\dz_1^{-1}\approx 0.212077,$$
 the iteration
arrives at the expected result at the second step $(n=2)$:
$$0.65-z_0\approx 0.437923;\qd 0.65-z_1\approx  0.430603;\qd 0.65-z_2\approx  0.430408.$$
}\dexmp

\prf First, we have $m=65/100$ and then
$$Q=\begin{pmatrix} -2/5 &  2/5\\
7/50 & -53/100\end{pmatrix}.$$
In this case, we ignore $(c_i)$ but let $b_1>0$. Actually, we have
$$\gathered
b_0= 2/5, \qd b_1= 39/100;\qd a_1= 7/50;\\
\mz_0= 1,\;  \mz_1 = 20/7; \qd \fz_0 = 265/78,\; \fz_1= 35/39.\endgathered$$
Therefore
$$v_0=\big(\sqrt{53/67},\; \sqrt{14/67\,}\big),\qqd
z_0^{-1}=5 \big(2809 + 40 \sqrt{742}\,\big)/4134.$$
The conclusion now follows by our algorithm.
\deprf

An additional example for the algorithm presented in this section is delayed
to Example \ref{4-21}.

Before moving further, let us introduce an algorithm for (and then
a representation of) the solution to equation (\ref{3-1}). This is mainly
used in theoretic analysis rather than numerical computation.
The idea is meaningful in a more general setup and comes from
\rf{chzh}{Theorem 1.1 and Proposition 2.6} plus
a modification \rf{cmf14}{Proposition 4.1}.
Given a number $z\in {\mathbb R}$ and a vector $v$, consider the equation for the vector $w$:
\be Qw +z w=-v.\lb{3-06}\de
To do so, we need some notation.
Fix $i: 0\le i\le N-1$, and set
$${\az}_{\ell}^{(i)}=\frac{1}{b_{i+\ell}}
\begin{cases}
c_{i+\ell}-z+ a_{i+\ell}, \qd & 1=\ell\le N-i\\
c_{i+\ell}-z, \qd & 2\le \ell \le N-i.
\end{cases}$$
Next, define the vector $G_{\cdot, 1}^{(i)}$ by $G_{\ell, 1}^{(i)}={\az}_{\ell}^{(i)}$ for $\ell=1, 2, \ldots, N-i$
and define recursively in $k=2,\;3,\; \ldots,\; N-i$, the vector $G_{\cdot,\, k}^{(i)}$ by
\begin{align}
&G_{\ell, k}^{(i)}= G_{\ell,\, k-1}^{(i)}+ {\az}_{\ell -k+1}^{(i+k-1)}G_{k-1,\, k-1}^{(i)},\qquad
\ell= k,\; k+1,\;\ldots,\; N-i. \lb{g-01}
\end{align}
Note that here for computing $G_{\cdot,\, k}^{(i)}$, we use only $G_{\cdot,\, k-1}^{(i)}$ but not the others $G_{\cdot,\, j}^{(i)}$ with $j\le k-2$.

\prp\lb{t-08-1}\;\;{\cms Let $N\ge 1$ and
 $G_{0,0}^{(\cdot)}\equiv 1$. Then the solution $w=(w_k: k\in E)$ to equation $(\ref{3-06})$ has the following representation:
$$w_n=\frac{v_N+M_{N-1}(v)}{c_N-z +M_{N-1}(c_{\cdot}-z)}[1+N_{n-1}(c_{\cdot}-z)]-N_{n-1}(v),\qqd 0\le n\le N$$
where for each vector $h$,  $N_{-1}(h)=0$ and
$$\begin{gathered}
M_{N-1}(h)= c_N\sum_{j=0}^{N-1}\frac{h_j}{b_{j}} G_{N-j,\, N-j}^{(j)},\\
N_n(h)= \sum_{j=0}^{n} \frac{h_j}{b_{j}}\sum_{k=0}^{n-j} G_{k,\,k}^{(j)},
\qqd 0\le n< N.\end{gathered}$$
}\deprp

The proof of this result is delayed to Section \ref{s-05}.

From now on, we are going to treat general real matrices.
This is a hard task and will be the main goal of the next
section. Here we study a special case only.
In computational mathematics, there is a well-known
Lanczos tridiagonalization procedure making a matrix
to be tridiagonal one. That is, for a given $A$,
constructing a nonsingular $B$ such that
$B^{-1} A B=:T$ becomes a tridiagonal matrix. We will come
back to the procedure soon. Here is an example
(the details are omitted).

\xmp \lb{t-09}\;\;{\cms Let
$${A=\left(\begin{matrix}1 & 2 & 3\\
1&2&1\\
3 & 2& 1\end{matrix}\right),}\qqd
{B=\left(\begin{matrix} 1 & 0 & 0\\
0&1/ \sqrt{10} & 3 /\sqrt{13}\\
0 & 3/ \sqrt{10}& -2/\sqrt{13}\end{matrix}\right).}$$
Then
$$\gathered
T=B^{-1}A B=\left(\begin{matrix}1 & 11/\sqrt{10} & 0\\
\sqrt{10}&25/11& 20\sqrt{130}/143 \\
0 & \sqrt{130}/11& 8/11\end{matrix}\right),\\
\rho (A)=\rho(T)=3 + \sqrt{5}\approx 5.23607.\endgathered$$
Our algorithm
arrives at the same result at the second step of the iterations $(n=2)$:
$$(n=0)\; 5.43937;\qd (n=1)\; 5.23996;\qd (n=2)\; 5.23607.$$
}\dexmp

It is the position to recommend an improved algorithm as follows.
The point is to use the inner product $(\cdot, \cdot)_{\mu}$ and norm $\|\cdot\|_{\mu}$
in the space $L^2(\mu)$ since $(\mu_k)$ may not be a constant as in Example \ref{t-01}.
\smallskip

\nnd{\bf Improved algorithm}.\; Given ${\tilde v}_0$ and $\dz_1$ as above,
redefine $v_0={\tilde v}_0/\|{\tilde v}_0\|_{\mu}$ and
$$z_0={\xi}{\dz_1^{-1}}+ (1-\xi) (v_0, -Q v_0 )_{\mu},\qqd \xi\in [0, 1].$$
For $k\ge 1$, define $w_k$ as before but redefine
$$v_k=\frac{w_k}{\|w_k\|_{\mu}}, \qqd z_k=(v_k, -Q v_k )_{\mu}.$$
\smallskip\

With $\xi=7/8$, Example \ref{t-07} and Table 2 are improved as Table 2'.

\begin{center}Table 2'.\; Example \ref{t-07} and Table 2 are improved using the new $z_0$ with $\xi\!=\!7/8$.\end{center}
\vskip-0.5truecm
$$
\hfil\vbox{\hbox{\vbox{\offinterlineskip
 \halign{&\vrule#&\strut\quad\hfil#\hfil\quad&\vrule#&
 \quad\hfil#\hfil\quad&\vrule#&
 \quad\hfil#\hfil\quad\cr
    \noalign{\hrule}
   height2pt&\omit&&\omit&&\omit&&\omit&&\omit&\cr
 & $\pmb{N+1}$ && $\pmb{z_0}$ && $\pmb{z_1}$ && {$\pmb{z_2=\lz_0}$} &&\pmb{\!\! upper/lower\!\!} &\cr
   height2pt&\omit&&\omit&&\omit&&\omit&&\omit&\cr
   \noalign{\hrule}
   height2pt&\omit&&\omit&&\omit&&\omit&&\omit&\cr
& $8$ && 0.523309 && {0.525268} && 0.525268 && $1\!+\!10^{-11}$&\cr
   height2pt&\omit&&\omit&&\omit&&\omit&&\omit&\cr
   \noalign{\hrule}
   height2pt&\omit&&\omit&&\omit&&\omit&&\omit&\cr
& $100$ && 0.387333 && 0.376393 && 0.376383 && $1\!+\!10^{-8}$  &\cr
  height2pt&\omit&&\omit&&\omit&&\omit&&\omit&\cr
   \noalign{\hrule}
   height2pt&\omit&&\omit&&\omit&&\omit&&\omit&\cr
& $500$ && 0.349147 && 0.338342 && 0.338329 && $1\!+\!10^{-7}$ &\cr
  height2pt&\omit&&\omit&&\omit&&\omit&&\omit&\cr
   \noalign{\hrule}
   height2pt&\omit&&\omit&&\omit&&\omit&&\omit&\cr
&\!$1000$\! && 0.338027 && 0.327254 && 0.32724 && $1\!+\!10^{-7}$ &\cr
  height2pt&\omit&&\omit&&\omit&&\omit&&\omit&\cr
   \noalign{\hrule}
   height2pt&\omit&&\omit&&\omit&&\omit&&\omit&\cr
& \!$5000$\! && 0.319895 && 0.30855 && 0.308529 && $1\!+\!10^{-7}$ &\cr
  height2pt&\omit&&\omit&&\omit&&\omit&&\omit&\cr
   \noalign{\hrule}
   height2pt&\omit&&\omit&&\omit&&\omit&&\omit&\cr
& \!$7500$\! && 0.316529 && 0.304942 && 0.304918 && $1\!+\!10^{-7}$ &\cr
  height2pt&\omit&&\omit&&\omit&&\omit&&\omit&\cr
   \noalign{\hrule}
   height2pt&\omit&&\omit&&\omit&&\omit&&\omit&\cr
& $10^4$ && 0.31437 && 0.302586 && 0.302561 && $1\!+\!10^{-7}$ &\cr
     height2pt&\omit&&\omit&&\omit&&\omit&&\omit&\cr
 \noalign{\hrule}}}}}\hfill$$
The last column is the order of the ratio of the upper and lower bounds of
$\lz_0$ in terms of $v_2$, as will be explained below, above Example \ref{t-12}.

Table 3 gives two more examples.

\begin{center}Table 3.\quad Outputs using improved  $z_0$ with $\xi=7/8$.\end{center}
\vskip-0.5truecm
$$
\hfil\vbox{\hbox{\vbox{\offinterlineskip
 \halign{&\vrule#&\strut\quad\hfil#\hfil\quad&\vrule#&
 \quad\hfil#\hfil\quad&\vrule#&
 \quad\hfil#\hfil\quad\cr
    \noalign{\hrule}
   height2pt&\omit&&\omit&&\omit&&\omit&\cr
   & \pmb{Example} && $\pmb {z_0}$ && $\pmb {z_1}$ &&  $\pmb {z_2=\lz_0}$ &\cr
   height2pt&\omit&&\omit&&\omit&&\omit&\cr
   \noalign{\hrule}
   height2pt&\omit&&\omit&&\omit&&\omit&\cr
 & \ref{t-08} && $0.436733$ && $0.430407$ && $0.430408$
    &\cr
   height2pt&\omit&&\omit&&\omit&&\omit&\cr
   \noalign{\hrule}
   height2pt&\omit&&\omit&&\omit&&\omit&\cr
   & \ref{t-09} && $5.36161$ && $5.23578$ && $5.23607$
    &\cr
   height2pt&\omit&&\omit&&\omit&&\omit&\cr
 \noalign{\hrule}}}}}\hfill$$

\subsection*{Appendix of \S\,\ref{s-03}. Algorithm for Lanczos tridiagonalization}
For a given $A$, the aim is choosing a nonsingular $Q$ such that
$$Q^{-1} A Q=T=
\left(\begin{matrix}
c_1 & b_1 & \cdots&\cdots & 0\\
a_1 & c_2 &  \ddots&{}   &\vdots \\
\vdots& \ddots & \ddots& \ddots &\vdots\\
\vdots & {} &\ddots &\ddots  & b_{n-1}\\
0& \cdots&\cdots & a_{n-1} & c_n
\end{matrix}\right).$$
Note that the notation here is somehow different from the other
part of the paper. To do so, we use the following column partitionings:
$$\aligned
Q&=[q_1|\cdots|q_n],\qqd
(Q^{-1})^*={\widetilde Q}=[{\tilde q}_1|\cdots|{\tilde q}_n].
\endaligned$$

Let $q_0=0$, ${\tilde q}_0=0$, $b_0=0$ and $a_0=0$.
Choose unit vectors $q_1$ and ${\tilde q}_1$ such that ${\tilde q}_1^* q_1=1$.
Define
$$\aligned
c_k&= {\tilde q}_k^* A q_k,\qqd k\ge 1,\\
r_k&=(A-c_k I) q_k-a_{k-1} q_{k-1},\qqd k\ge 1,\\
{\tilde r}_k&=(A-c_k I)^* {\tilde q}_k-b_{k-1} {\tilde q}_{k-1},\qqd k\ge 1,\\
b_k&=\|r_k\|_2, \qqd k\ge 1,\\
a_k&= {\tilde r}_k^* r_k/b_k, \qqd k\ge 1,\\
q_k&=r_{k-1}/b_{k-1},\qqd k\ge 2,\\
{\tilde q}_k&={\tilde r}_{k-1}/a_{k-1},\qqd k\ge 2.
\endaligned
$$

For Example \ref{t-09}, we simply choose $q=(1,0,0)^*$ and ${\tilde q}=(1,0,0)^*$. Gene\-rally speaking, there is a question in choosing initial
$q_0$ and ${\tilde q}_0$. More generally, it should be meaningful to know
for what $A$, the resulting matrix have positive $a_k$ and $b_k$ for every $k$.

\section{Efficient initials. General case}\lb{s-04}

A general algorithm for the efficient initials will be introduced
later in the second subsection.
The algorithm introduced in the next subsection is easier and quite
general, but may be less efficient.

\subsection{Fix uniformly distributed initial vector {$\pmb {v_0}$}}\lb{s-4-1}

In this subsection, we fix the uniformly distributed initial vector
$$v_0=(1, 1, \ldots, 1)/\sqrt{N+1}.$$
This is the easiest choice of $v_0$ since it does not use any information
from the eigenvector $g$ of $\rho (A)$ except its positivity property.
On the other hand, this means that the choice is less efficient and it can
be even broken as shown by Example \ref{t-03}. The effectiveness of this $v_0$
depends heavily on the choice of $z_0$. For which, here we introduce
three effective choices.
\medskip

\nnd{\bf Choice I}\qd Let $A\!=\!(a_{ij}: i,j\!\in\! E)$ be nonnegative and set $z_0\!=\!\sup_{i\in E} A_i$, where $A_i=\sum_{j\in E}a_{ij}$. This universal choice comes from the
fact that $\sup_{i\in E} A_i$ is a upper bound of $\rho (A)$, which
can be seen by setting $x_i\equiv 1$ in the next result.

\prp\lb{t-11}\;\;{\cms For a nonnegative irreducible matrix $A$ with maximal eigenvalue $\rho (A)$,
the Collatz--Wielandt formula holds:
$$\sup_{x>0}\min_{i\in E} \frac{(Ax)_i}{x_i}=\rho (A)
= \inf_{x>0}\max_{i\in E} \frac{(Ax)_i}{x_i}.$$
}\deprp

For the present $(v_0, z_0)$, even though it is not necessary, one may replace (\ref{02}) by
\be (z_{k-1}I-A) w_k =v_{k-1}. \de
This choice of $z_0$ avoids the  collapse of the algorithm since
$$0<z_0-\rho(A)< |z_0 -\lz|$$
for every eigenvalue $\lz\ne \rho(A)$ of $A$.

Let us now introduce an important application of Proposition \ref{t-11}.
First, if we replace $A$ and $\rho(A)$ with $-Q$ and $\lz_0$, respectively, the same conclusion holds, as shown in the next corollary (the proof is delayed to Section \ref{s-05}). Actually, the corollary holds in a much more general setup. Refer to
\rf{cmf04}{Theorem 9.5}.

\crl\lb{t-11-1}\;\; {\cms For $Q$-matrix, the Collatz--Wielandt formula becomes
$$ \sup_{x>0}\min_{i\in E} \frac{(-Qx)_i}{x_i}=\lz_0 (Q)
= \inf_{x>0}\max_{i\in E} \frac{(-Qx)_i}{x_i}.$$
}\decrl

Thus, instead
of the mean estimate given in these algorithm, we can produce pointwise
estimates. To do so, we need only to compute the ratio
$(-Q)v_k/v_k$. For instance, in Example \ref{t-02}, the ratio $(-Q)v_2/v_2$ is as follows:
$$0.525197,\, 0.5254,\, 0.52553,\, 0.525623,\, 0.525693,\, 0.525747, \, 0.525787,\,0.525816.$$
Therefore, we obtain
$$0.525197\le \lz_0\le 0.525816$$
and the ratio of the upper/lower bounds is $\approx 1.00118$.
Next, for Example \ref{t-07},
the ratio $(-Q)v_2/v_2$ is as follows:
$$0.525268,\,0.525268,\,0.525267,0.525267,0.525267,0.525267,0.525267,0.525267.$$
Hence, we have
$$0.525267\le \lz_0\le 0.525268$$
and the ratio of the upper/lower bounds is $\approx 1+  10^{-6}$.
Actually, if we apply the estimates given in \rf{cmf10}{Theorem 2.4\,(3)}
(with supp\,$(f)=E$):
$$\aligned
&z_2\wedge \Big(\sup_{i\in E} f_i/g_i\Big)\ge \lz_0 \ge \inf_{i\in E} f_i/g_i,\\
&g_i:= \sum_{k\in E} \mu_k f_k \fz_{i\vee k}
=\fz_i\sum_{k=0}^i \mu_k f_k + \sum_{i+1\le k\le N} \mu_k \fz_k f_k,\\
&\fz_i:=\sum_{k=i}^N \frac{1}{\mu_k b_k}\qqd (\text{for this example }\mu_k\equiv 1,\; b_i=(i+1)^2),
\endaligned$$
to the test function $f=v_2$ with a more precise output, the upper/lower bounds can be improved as $\approx 1+ 10^{-7}$. Hence the estimate $\lz_0\approx 0.525268$ is indeed sharp up to
the six precisely significant digits. This shows that
the estimates in the latter example are much better than the former one.

\xmp \lb{t-12}\;\;{\cms Let $A$ be the same as in Example $\ref{t-09}$. Then
$\rho (A)\approx 5.23607$ and $z_0=6$.
The Rayleigh quotient iteration gives us
$$z_1\approx 5.27273, \; z_2\approx 5.23639,\; z_3\approx 5.23607.$$
}\dexmp

\xmp\lb{t-13}\;\;{\cms Let
$$A=
\begin{pmatrix}
1 & 2 & 3 & 4\\
5 & 6 & 7 & 8\\
9 & 10 & 11 & 12\\
13 & 14 & 15 & 16
\end{pmatrix}.$$
Then
$\rho (A)\approx 36.2094$ and $z_0=58$.
The Rayleigh quotient iteration gives us
$$z_1\approx 37.3442, \; z_2\approx 36.2674,\; z_3\approx 36.2095,\;
z_4\approx 36.2094.$$
}\dexmp

\xmp\lb{t-14}\;\;{\cms Let
$$A=
\begin{pmatrix}
1 & 2 & 0 & 0\\
3 & 14 & 11 & 0\\
9 & 10 & 11 & 1\\
5 & 6 & 7 & 8
\end{pmatrix}.$$
This matrix has complex eigenvalues:
$$24.0293,\; 7.72254,\; 1.1241 + 2.40522\, i,\; 1.1241 - 2.40522\, i.$$
Hence
$\rho (A)\approx 24.0293$ and $z_0=31$.
The Rayleigh quotient iteration gives us
$$z_1\approx 24.4393, \; z_2\approx 24.0385,\; z_3\approx 24.0293.$$
}\dexmp

\xmp \lb{t-16}\;\;{\cms Let $Q$ be the same as in Example $\ref{t-01}$ and let
$$A=113\, I +Q.$$
Then $z_0=113$. Recall that $\lz_{\min} (-Q)\approx 0.525268$.
For $k=1, 2, 3$, the Rayleigh quotient iteration gives us $113-z_k$ as follows
$$113-z_1\approx 0.602312, \; 113-z_2\approx 0.525463,\; 113-z_3\approx 0.525268.$$
Alternatively, one may apply the algorithm directly to $-Q$ with $z_0=0$.}\dexmp

We remark that the algorithm is meaningful for any
$z_0\ge \sup_{i\in E}\sum_{j\in E}a_{ij}$. For instance, if we choose $z_0=200$         
rather than $z_0=6$ used in Example \ref{t-12}, then the successive results of the
iterations are as follows:
$$z_1\approx 5.33546,\; z_2\approx 5.24182,\; z_3\approx 5.23608,\; z_4\approx 5.23607.$$
The convergence becomes slower as we can imagine. In other words, a larger initial $z_0$ is less efficient. In view of Proposition \ref{t-11},
we have $0<\rho(A)\le 113$. It seems there is a large room for us to
choose $z_0$. Yes or no? It is yes, since the last estimates are rather rough,
each choice $z_0\in [111.7, 113]$ is also available.
The answer is also no, since if we choose $z_0=111.6$, then we will go to the
pitfall $\lz_1(>\lz_0)$. Hence, it is rather sensitive to find a
useful $z_0$ except the Choice I. Noting that $\rho(A)\approx 113-0.525268$, the reason why the rough Choice I is still efficient for this model should be clear.

We have thus studied the model introduced in Example 1 six times with
different initials. The results are collected in Table 4. Among them, the worst one is Example \ref{t-03} and the best one is Example \ref{t-07} which
uses the whole power of the algorithm introduced in Section \ref{s-03}.
The ``Uniform'' is the present Choice I and the ``Auto'' means automatic one
given by the algorithm, as we will come back in Choice II below.

\begin{center}Table 4.\quad Comparison of examples with different initials.\end{center}
\vspace{-1.1truecm}

$$
\hfil\vbox{\hbox{\vbox{\offinterlineskip
 \halign{&\vrule#&\strut\quad\hfil#\hfil\quad&\vrule#&
 \quad\hfil#\hfil\quad&\vrule#&
 \quad\hfil#\hfil\quad\cr
    \noalign{\hrule}
   height2pt&\omit&&\omit&&\omit&&\omit&\cr
   & {\pmb{Same $Q$}} && ${\pmb{v_0}}$ && ${\pmb{z_0}}$ &&  {\pmb{\# of Iterations}} &\cr
   height2pt&\omit&&\omit&&\omit&&\omit&\cr
   \noalign{\hrule}
   height2pt&\omit&&\omit&&\omit&&\omit&\cr
   & Example \ref{t-01} && ${\tilde v}_0$ && Power && $10^3$  &\cr
   height2pt&\omit&&\omit&&\omit&&\omit&\cr
   \noalign{\hrule}
   height2pt&\omit&&\omit&&\omit&&\omit&\cr
& Example \ref{t-02} && ${\tilde v}_0$ && Auto && $2$   &\cr
   height2pt&\omit&&\omit&&\omit&&\omit&\cr
   \noalign{\hrule}
   height2pt&\omit&&\omit&&\omit&&\omit&\cr
&Example \ref{t-03} && Uniform && Auto && Collapse   &\cr
   height2pt&\omit&&\omit&&\omit&&\omit&\cr
 \noalign{\hrule}
height2pt&\omit&&\omit&&\omit&&\omit&\cr
&Example \ref{t-04} && Uniform && $\dz_1^{-1}$ && $2$   &\cr
   height2pt&\omit&&\omit&&\omit&&\omit&\cr
 \noalign{\hrule}
 height2pt&\omit&&\omit&&\omit&&\omit&\cr
&Example \ref{t-07} && ${\tilde v}_0$ &&$\dz_1^{-1}$ && $2$   &\cr
   height2pt&\omit&&\omit&&\omit&&\omit&\cr
 \noalign{\hrule}
 height2pt&\omit&&\omit&&\omit&&\omit&\cr
&Example \ref{t-16} && Uniform && $113$ && $3$   &\cr
   height2pt&\omit&&\omit&&\omit&&\omit&\cr
 \noalign{\hrule}
}}}}\hfill$$

In conclusion, even though the present choice $(v_0, z_0)$ may not be very
efficient, but it works in a very general setup. This algorithm works even
for a more general class of matrices, without assuming the nonnegative property,
once you have a upper estimate of the largest eigenvalue of $A$.
Clearly, for large-scale matrix, Choice I is meaningful only for
the sparse ones.
\medskip

\nnd{\bf Choice II}\;\;Simply use the particular choice given in the
Rayleigh quotient iteration: $z_0=v_0^*Av_0$. This simple choice is
quite natural and so is often used in practice. However, there is a dangerous
here since $v_0$ is chosen roughly, the algorithm may lead to an incorrect limit, as illustrated by Example \ref{t-03}.

With the present $z_0$, the computation results for Examples \ref{t-12}--\ref{t-14} are listed in Table 5.
\vskip-0.5truecm

\begin{center}Table 5.\quad Output $(z_1, z_2, z_3)$ of Examples \ref{t-12}--\ref{t-14}.\end{center}
\vskip-0.5truecm
$$
\hfil\vbox{\hbox{\vbox{\offinterlineskip
 \halign{&\vrule#&\strut\quad\hfil#\hfil\quad&\vrule#&
 \quad\hfil#\hfil\quad&\vrule#&
 \quad\hfil#\hfil\quad\cr
    \noalign{\hrule}
   height2pt&\omit&&\omit&&\omit&&\omit&\cr
   & \pmb{Example} && $\pmb {z_1}$ && $\pmb {z_2}$ &&  $\pmb{z_3=\lz_0}$ &\cr
   height2pt&\omit&&\omit&&\omit&&\omit&\cr
   \noalign{\hrule}
   height2pt&\omit&&\omit&&\omit&&\omit&\cr
   & \ref{t-12} && $5.24183$ && $ 5.23608$ && $5.23607$
  &\cr
   height2pt&\omit&&\omit&&\omit&&\omit&\cr
   \noalign{\hrule}
   height2pt&\omit&&\omit&&\omit&&\omit&\cr
& \ref{t-13} && $35.8428$ && $36.2127$ && $36.2094$
    &\cr
   height2pt&\omit&&\omit&&\omit&&\omit&\cr
   \noalign{\hrule}
   height2pt&\omit&&\omit&&\omit&&\omit&\cr
   & \ref{t-14} && $23.7316$ && $24.0317$ && $24.0293$
    &\cr
   height2pt&\omit&&\omit&&\omit&&\omit&\cr
 \noalign{\hrule}}}}}\hfill$$
Combining $(z_1, z_2)$ here with those given in the last part,
it is clear that the present choice of $z_0$, once works,
 is better than Choice I.
\medskip

\nnd{\bf Choice III}\;\; This is based on a comparison technique.
For given $A=(a_{ij})$ having the property $a_{i, i+1}+a_{i+1, i}>0$
for every $i$, we introduce the symmetrized matrix $(A+A^*)/2.$
(This symmetrizing procedure may be omitted if both $a_{i, i+1}>0$ and $a_{i+1, i}>0$ for every $i$).
Denote by $(\az_i, \bz_i, \gz_i)$ the tridiagonal part
(where $\gz_i$ are the diagonal elements) taken from
the symmetrized matrix. By assumption, we have $\az_i>0$ and $\bz_i>0$.
We can then follow the last section to choose a $z_0$ first for
the tridiagonal matrix and then regarding it as an approximation of
$z_0$ for the original $A$. One may worry that we have lost too much in the last
step. Yes, it may be so. However, the key is to avoid the collapse. The
smaller estimate $z_0$ of $\lz_{\min}(-Q)$ is not really serious since the
algorithm can repair it rapidly, as shown by the next example.

\xmp \lb{t-17}\;\;{\cms Let $A$ be the same as in Example $\ref{t-14}$.
Then
$$\frac 1 2 (A+A^*)=\begin{pmatrix}
1 & 5/2 & 9/2 & 5/2\\
5/2 & 14 & 21/2 & 3\\
9/2 & 21/2 & 11 & 4\\
5/2 & 3 & 4 & 8
\end{pmatrix}.$$
From this, we obtain a tridiagonal matrix
$$T=\begin{pmatrix}
1 & 5/2 & 0 & 0\\
5/2 & 14 & 21/2 & 0\\
0 & 21/2 & 11 & 4\\
0 & 0 & 4 & 8
\end{pmatrix}$$
and then
$$Q=T-27\,I=\begin{pmatrix}
-26 & 5/2 & 0 & 0\\
5/2 & -13 & 21/2 & 0\\
0 & 21/2 & -16 & 4\\
0 & 0 & 4 & -19
\end{pmatrix}.$$
According to what we did in \S \ref{s-03}, we have
$z_0\approx 1/0.321526$ for $-Q$. Then, we have $z_0\approx 27- 1/0.321526$
for $T$. This is regarded as an approximation of $z_0$ for $A$. Starting from
here and using the Rayleigh quotient iteration, we obtain the successive
approximation of $\rho (A)$ as follows:
$$z_1\approx 24.0125, \quad z_2\approx 24.0293$$
as we expected. Picking up the tridiagonal part directly from $A$
(without using the symmetrizing procedure), the same approach leads to the following output:
$$z_0\approx 28-1/0.23307\approx 23.7094,\qd z_1\approx 23.9901,\qd z_2\approx 24.0293. $$}\dexmp

Let us remark that the three choices of $z_0$ in this subsection are independent of
the initial $v_0$ used here and so can be also used in the next subsection. Certainly, there are other approaches can be used to
deduce an approximation of the required $z_0$. For instance, the Cheeger's
approach \rf{cmf04}{\S 9.5}, which is mea\-ningful in a very general
setup. Since it takes account of all subset of $E$ (except the emptyset), the number of computations is of order $2^N$. This approach as well as
the capacitary one (cf. \rf{cmf05}{Chapter 7}) needs to be simplified to fit the present setup. In practice, one often
uses Proposition \ref{t-11} or Corollary \ref{t-11-1} to get a upper/lower
bound in terms of a suitable test sequence $(x_i)$. Refer also to \rf{cmf05}{Theorem 3.6} which uses test weights. These approaches depend heavily on the working models.

\subsection{Efficient initial vector {$\pmb {v_0}$}}\lb{s-4-2}

\medskip

In general, it is much more difficult to choose an efficient initial $v_0$ than $z_0$. Here is our algorithm.

\subsection*{A general algorithm}

Let $A=(a_{ij}: i, j\in E)$ be a given irreducible matrix
having nonnegative off-diagonal elements. Once again, denote by $\rho(A)$ the maximal eigenvalue of $A$. If $A_i:=\sum_{j\in E} a_{ij}$ is a constant (independent of $i\in E$), then we have $\rho(A)\equiv A_i$ with right-eigenvector $\bbb {1}$ (its components are all equal to $1$). From now on, we assume that $A_i$ are not a constant.

We introduce our algorithm is four steps.

{\it Step} 1. If $A_i\le 0$ for every $i\in E$, one can jump from here to Step 2 below
by setting $Q=A$. Otherwise, let $\max_{i\in E} A_i>0$. Define
$$Q= A- \Big(\max_{i\in E} A_i\Big)I.$$
Then the sum of each row of $Q$ is less or equal to zero and at least one of the rows is less than zero since $A_i$ is not a constant. Now, if
$$Q_0=\cdots =Q_{N-1}=0$$ but
$Q_N\!<\!0$ ($Q_k\!:=\!\sum_{j}q_{kj}$), then one can jump from here to Step 3 with $h_i\equiv 1$.

{\it Step} 2. Assume that $Q_k<0$ for some $k\le N-1$. Denote by
$h=(h_0, h_1,\ldots, h_N)^*$ with $h_0=1$  the solution to the equation
$$Q^{\setminus N\text{'s row}}h=0,$$
where $Q^{\setminus \text{$k$'s row}}$ is obtained from $Q$ removing its $k$'s row $(q_{k0}, q_{k1},\ldots, q_{kN})$.
In the case that
$$c_N+\sum_{j\le N-1}q_{Nj}\bigg(1-\frac{h_j}{h_N}\bigg)$$
is much smaller then
$$\sum_{j\le N-1}q_{Nj} \frac{h_j}{h_N}$$
(say, $1:100$ for instance), one can jump from here to (\ref{4-09}) with $x_i\equiv 1$ (cf. Example \ref{4-20} in the case of $b_4=0.01$).

{\it Step} 3. Let $(h_i: i\in E)$ be constructed in  the last step.
Define $q_i=-q_{ii}$, $i\in E$.
Let $x=(x_0, x_1,\ldots, x_N)^*$ (with $x_0=1)$ be the solution to
the equation
\be x^{\setminus \text{$0$'s row}}
=P^{\setminus \text{$0$'s row}}\,x,    \lb{4-8-1}\de
where
$$P=(p_{ij}: i,j\in E): \qqd p_{ii}=0,\qd p_{ij}=\frac{q_{ij}h_j}{q_i h_i}, \qqd j\ne i;$$
or in the matrix form,
$$P=\text{\rm Diag}\big((q_i h_i)^{-1}\big)Q\,\text{\rm Diag}(h_i)+I.$$
Refer to the comments below Examples \ref{4-20} and \ref{4-21} for the constraint
$x_0=1$.
Here the sequence $(x_i)$ is an extension
of $(\fz_i)$ used in Section \ref{s-03} (cf. Lemma \ref{t-21} below).

{\it Step} 4. We are now ready to state our algorithm as follows.
Define a (column) vector ${\tilde v}_0$ with components
\be {\tilde v}_0(i)=h_i \sqrt{x_i}, \qqd i=0,1,\ldots, N. \lb{4-09} \de
Let
$$\aligned
v_0&=\frac{{\tilde v}_0}{\sqrt{{\tilde v}_0^*{\tilde v}_0}},\qquad
z_0=v_0^* (-Q)v_0.
\endaligned
$$
In general, for $k\ge 1$, let $w_{k}$ be the solution to the
equation
$$(-Q-z_{k-1} I)w_{k}=v_{k-1},$$
and define
$$\aligned
v_{k}&=\frac{w_{k}}{\sqrt{w_{k}^*w_{k}}},\qquad
z_{k}=v_{k}^* (-Q)v_{k}.
\endaligned
$$
Then $z_{k}$ and $v_k$ are approximations of the minimal eigenvalue
 $\lz_0=\lambda_{\min}(-Q)$ of $-Q$ and its eigenvector, respectively. If we replace $-Q$ by $A$ everywhere in this step, then the resulting $z_{k}$ and $v_k$ are approximations of $\rho(A)$ and its eigenvector $g$, respectively. Obviously, from Step 1, it
 follows that
$$\lambda_{\min}(-Q)+\rho(A)=\max_{i\in E}A_i.$$
Hence,
$$\lz_0=\lambda_{\min}(-Q)>\az \Longleftrightarrow \rho(A)
\le \max_{i\in E}A_i-\az.$$
This gives the relationship of a lower estimate of $\lz_0$ and
an upper estimate of $\rho(A)$.

\xmp\lb{x-18} \;\;{\cms Let $A$ be given in Example $\ref{t-09}$. Then
$$\rho(A)=3 + \sqrt{5}\approx 5.23607.$$
Our algorithm here gives us
$$z_1\approx 5.23883, \qd z_2\approx  5.23607.$$}
\dexmp

\prf Since $\max_i A_i=6$, we have
$$Q=A- 6 \,I
= \begin{pmatrix}
-5 & 2 & 3\\
1 & -4 & 1\\
3 & 2 & -5
\end{pmatrix}.$$
Next, we have
$$h_0=1,\qd h_1= 4/7,\qd h_2= 9/7$$
and
$$x_0=1,\qd x_1= 7/9,\qd x_2 = 49/81.$$
From these, we obtain
$${\tilde v}_0=\big(1,\; h_1\sqrt{x_1},\; h_2\sqrt{x_2}\,\big)^*
=\big(1,\; 4/\big(3 \sqrt{7}\,\big),\; 1 \big)^*.$$
Now, with
$$v_0={\tilde v}_0/\sqrt{{\tilde v}_0^*{\tilde v}_0},\qqd
z_0=v_0^*A v_0\approx 5.11616,$$
we can apply the Rayleigh quotient iteration
in two steps to obtain the conclusion.\deprf

\xmp\lb{x-19} \;\;{\cms Let $A$ be the same as in Example $\ref{t-13}$.
Then $\rho(A)\approx 36.2094$.
By using $(\ref{4-09})$,
$$v_0=(0.348213,\; 0.244601,\; 0.389728,\; 0.816719)^*,$$
the Rayleigh quotient iteration starts at $z_0\approx 34.4924$
and gives us
$$z_1\approx 36.1469,\;\; z_2\approx 36.2095,\;\; z_3\approx 36.2094.$$
}\dexmp

\prf We have
$$Q= A - 58\,I=
\begin{pmatrix}
-57 & 2 & 3 & 4 \\
5 & -52 & 7 & 8\\
9 & 10 & -47 & 12\\
13 & 14 & 15 & -42
\end{pmatrix}.$$
Next, we have
$$h_0=1,\;\; h_1= \frac{59}{27},\;\; h_2 = \frac{91}{27}, \;\;
h_3 = \frac{287}{27}.$$
Furthermore, we have
$$x_0=1,\;\; x_1 = \frac{189}{1829},\;\;
 x_2 = \frac{7155}{64883},\;\; x_3 = \frac{243}{4991}.$$
 Then the conclusion follows from the iteration.
\deprf

\xmp\lb{x-20} \;\;{\cms Let $A$ be the same as in Example $\ref{t-14}$.
Then $\rho(A)\approx 24.0293$.
By using the algorithm in Section $\ref{s-4-2}$,
the Rayleigh quotient iteration starts at $31-z_0\approx 13.7532$
and gives us for $k=1, \ldots, 4$,
$$31-z_k\approx 7.10985, \qd 7.60885, \qd 7.72193, \qd 7.72254, $$
respectively.}\dexmp

\prf We have
$$Q= A - 31\,I
=\begin{pmatrix}
-30 & 2 & 0 & 0\\
3 & -17 & 11 & 0\\
9 & 10 & -20 & 1\\
5 & 6 & 7 & -23
\end{pmatrix}.$$
Then, we have
$$\aligned
&h_0=1, \qd  h_1 = 15,\qd h_2 = 252/11,\qd  h_3 = 3291/11;\\
&x_0=1, \qd x_1 = 3691/76575,\qd x_2 = 1694/45945,\qd x_3 = 7447/3360111;\\
&v_0=(0.0659989,\; 0.217349,\; 0.290324,\; 0.929578).
\endaligned$$
The conclusion follows by the algorithm.\footnote{The author thanks
Xiao-Jun Mao who pointed out an error of $v_0$ in the published version of the paper.}\deprf

Due to Choice II of $z_0$, this example is in a similar situation as in Example \ref{t-03}. We will return to this example in Tables 11 and 12 below.

Actually, to show our algorithm is reasonable, one may ignore the
part using the $H$-transform and jump to the last step on $Q$-matrix
since the transform does not change the spectrum. Thus, one needs to
compare the maximal eigenvector $g$ and its approximation $(x_i)$. As
mentioned before, this depends heavily on the rate $b_N=c_N$. Here is
an example of sparse matrix.

\def\u01{\begin{matrix}(1,\; 1.00011,\; 1.00017,\; 0.999498,\; 0.998616)^*\\
(1,\; 1,\; 1,\; 0.999728,\; 0.999274)^*
\end{matrix}}
\def\v-1{\begin{matrix}(1,\; 1.00992,\; 1.0149,\; 0.955637,\; 0.877794)^*\\
(1,\; 1,\; 1,\; 0.9759,\; 0.934353)^*
\end{matrix}}
\def\x100{\begin{matrix}(1,\; 1.08011,\; 1.1211,\; 0.656961,\; 0.0652116)^*\\
(1,\; 1,\; 1,\; 0.805682,\; 0.253629)^*
\end{matrix}}
\def\y300{\begin{matrix}(1,\; 1.08399,\; 1.12701,\; 0.641158,\; 0.0226915)^*\\
(1,\; 1,\; 1,\;0.795865,\; 0.149677)^*
\end{matrix}}
\def\z500{\begin{matrix}(1,\; 1.08481,\; 1.12826,\; 0.637827,\; 0.0137353)^*\\
(1,\; 1,\; 1,\;0.79378,\; 0.116462)^*
\end{matrix}}

\xmp\lb{4-20}\;\;{\cms Let
$$Q=\begin{pmatrix}
-3 & 2 & 0 & 1 & 0\\
4 &  -7 &  3 & 0 & 0\\
0 & 5 & -5 & 0 & 0\\
10 & 0 & 0 & -16 &  6\\
0 & 0 & 0 & 11 & -11-b_4
\end{pmatrix}.  $$
Corresponding to different $b_4$, the
maximal eigenvector $g$ (normalized so that the first component to be one) and its approximation $(\sqrt{x_i}\,)$
(up to a positive constant) are given in Table $6$.
\vskip-0.5truecm

\begin{center}Table $6$.\quad  For different $b_4$, the
vectors $g$ and $(\sqrt{x_i}\,)$.
\end{center}

\vskip-0.6truecm

$$
\hfil\vbox{\hbox{\vbox{\offinterlineskip
 \halign{&\vrule#&\strut\quad\hfil#\hfil\quad&\vrule#&
 \quad\hfil#\hfil\quad&\vrule#&
 \quad\hfil#\hfil\quad\cr
    \noalign{\hrule}
   height2pt&\omit&&\omit&\cr
   & {$\pmb{b_4}$} && {${\pmb g}\,\text{\cms \big(1$^{\text{\cms st}}$ line\big)}$ and ${\pmb {\sqrt{x}}}\,\text{\cms \big(2$^{\text{\cms ed}}$ line\big)}$}\quad up to a constant &\cr
   height2pt&\omit&&\omit&\cr
   \noalign{\hrule}
   height2pt&\omit&&\omit&\cr
   & 0.01 && {$\displaystyle{\u01}$}
  &\cr
   height2pt&\omit&&\omit&\cr
   \noalign{\hrule}
   height2pt&\omit&&\omit&\cr
& 1 && {$\displaystyle{\v-1}$}
&\cr
   height2pt&\omit&&\omit&\cr
   \noalign{\hrule}
   height2pt&\omit&&\omit&\cr
& 100 && {$\displaystyle{\x100}$}
&\cr
   height2pt&\omit&&\omit&\cr
 \noalign{\hrule}}}}}\hfill$$
The corresponding output of our algorithm is given in Table $7$.
\vskip-0.5truecm

\begin{center}Table $7$.\quad  For different $b_4$, the
eigenvalue $\lz_0$ and $z_1, z_2, z_3$.
\end{center}
\vskip-0.4truecm

$$
\hfil\vbox{\hbox{\vbox{\offinterlineskip
 \halign{&\vrule#&\strut\quad\hfil#\hfil\quad&\vrule#&
 \quad\hfil#\hfil\quad&\vrule#&
 \quad\hfil#\hfil\quad\cr
    \noalign{\hrule}
   height2pt&\omit&&\omit&&\omit&&\omit&&\omit&\cr
   & {$\pmb{b_4}$}
   &&${\pmb{z_1}}$ &&${\pmb{z_2}}$ &&${\pmb{z_3=\lz_0}}$ &\cr
   height2pt&\omit&&\omit&&\omit&&\omit&&\omit&\cr
   \noalign{\hrule}
   height2pt&\omit&&\omit&&\omit&&\omit&&\omit&\cr
& \!0.01\!  &&\!0.000278573\!  && \!0.000278686\! && {}  &\cr
   height2pt&\omit&&\omit&&\omit&&\omit&&\omit&\cr
   \noalign{\hrule}
   height2pt&\omit&&\omit&&\omit&&\omit&&\omit&\cr
& 1  &&0.0236258  && 0.0245174 && 0.0245175  &\cr
   height2pt&\omit&&\omit&&\omit&&\omit&&\omit&\cr
   \noalign{\hrule}
   height2pt&\omit&&\omit&&\omit&&\omit&&\omit&\cr
& 100  &&0.200058  && 0.182609 && 0.182819  &\cr
    height2pt&\omit&&\omit&&\omit&&\omit&&\omit&\cr
 \noalign{\hrule}}}}}\hfill$$
}
\dexmp

Our original purpose to design the $Q$-matrix in the last example is
for a test of sparse matrix. The solution $x_0=x_1=x_2=1$ leads us to
think about the transition machinery of the $Q$-matrix.
Here is the graphic structure of the $Q$-matrix:
$$\cirnm{2}\leftrightarrows\cirnm{1}\leftrightarrows\cirnm{0}
\leftrightarrows\cirnm{3}\leftrightarrows\cirnm{4}.$$
As we will see at the end of Section \ref{s-05}, $x_i$ is the probability of the process first hitting $0$ starting from $i$ (which is exactly the probabilistic meaning of the construction of $(x_i)$ given in our general algorithm). Now,
starting from $2$, there is only one way to go to $0$, and hence
$x_2$ should be equal to $1$. So does $x_1$. From this graph,
it follows that the matrix is indeed tridiagonal after a relabeling
(simply exchange the labels $\cirnm{2}$ and $\cirnm{0}$):
$$\cirnm{0}\leftrightarrows\cirnm{1}\leftrightarrows\cirnm{2}
\leftrightarrows\cirnm{3}\leftrightarrows\cirnm{4}.$$
As a comparison, we present the next result using the
algorithms given in Sections \ref{s-04} and \ref{s-03}, respectively.

\def\uu01{\begin{matrix}(1,\; 0.999944,\; 0.999833,\; 0.999331, \; 0.998449)^*\\
(1,\; 0.999819,\; 0.999682,\; 0.99941,\; 0.998956)^*
\end{matrix}}
\def\vv-1{\begin{matrix}(1,\; 0.995096,\; 0.98532,\; 0.941608,\; 0.864908)^*\\
(1,\; 0.984848,\; 0.973329,\; 0.949871,\; 0.909433)^*
\end{matrix}}
\def\xx100{\begin{matrix}(1,\; 0.963436,\; 0.89198,\; 0.585996,\; 0.0581675)^*\\
(1,\; 0.91325,\; 0.842344,\; 0.678661,\; 0.213643)^*
\end{matrix}}

\xmp \lb{4-21}\;\;{\cms Let
$$Q=\begin{pmatrix}
-5 & 5 & 0 & 0 & 0\\
3 &  -7 &  4 & 0 & 0\\
0 & 2 & -3 & 1 & 0\\
0 & 0 & 10 & -16 &  6\\
0 & 0 & 0 & 11 & -11-b_4
\end{pmatrix}.  $$
Corresponding to different $b_4$, the
maximal eigenvector $g$ and its approximation $(\sqrt{x_i}\,)$
 are given in Table $8$.
\vskip-0.5truecm

\begin{center}Table $8$.\quad  For different $b_4$, the
vectors $g$ and $(\sqrt{x_i}\,)$.
\end{center}
\vskip-0.8truecm

$$
\hfil\vbox{\hbox{\vbox{\offinterlineskip
 \halign{&\vrule#&\strut\quad\hfil#\hfil\quad&\vrule#&
 \quad\hfil#\hfil\quad&\vrule#&
 \quad\hfil#\hfil\quad\cr
    \noalign{\hrule}
   height2pt&\omit&&\omit&\cr
   & {$\pmb{b_4}$} && {${\pmb g}\,\text{\cms \big(1$^{\text{\cms st}}$ line\big)}$ and ${\pmb {\sqrt{x}}}\,\text{\cms \big(2$^{\text{\cms ed}}$ line\big)}$}\quad up to a constant &\cr
   height2pt&\omit&&\omit&\cr
   \noalign{\hrule}
   height2pt&\omit&&\omit&\cr
   & 0.01 && {$\displaystyle{\uu01}$}
  &\cr
   height2pt&\omit&&\omit&\cr
   \noalign{\hrule}
   height2pt&\omit&&\omit&\cr
& 1 && {$\displaystyle{\vv-1}$}
&\cr
   height2pt&\omit&&\omit&\cr
   \noalign{\hrule}
   height2pt&\omit&&\omit&\cr
& 100 && {$\displaystyle{\xx100}$}
&\cr
   height2pt&\omit&&\omit&\cr
 \noalign{\hrule}}}}}\hfill$$
The corresponding output $(z_k)$ of the algorithm in Section $\ref{s-04}$ is given in Table $9$.
\vskip-0.5truecm

 \begin{center}Table $9$.\quad  For different $b_4$, the
eigenvalue $\lz_0$ and $z_1, z_2, z_3$.
\end{center}
\vskip-0.4truecm

$$
\hfil\vbox{\hbox{\vbox{\offinterlineskip
 \halign{&\vrule#&\strut\quad\hfil#\hfil\quad&\vrule#&
 \quad\hfil#\hfil\quad&\vrule#&
 \quad\hfil#\hfil\quad\cr
    \noalign{\hrule}
  height2pt&\omit&&\omit&&\omit&&\omit&&\omit&\cr
   & {$\pmb{b_4}$}
   &&${\pmb{z_1}}$ &&${\pmb{z_2}}$ &&${\pmb{z_3=\lz_0}}$ &\cr
   height2pt&\omit&&\omit&&\omit&&\omit&&\omit&\cr
   \noalign{\hrule}
  height2pt&\omit&&\omit&&\omit&&\omit&&\omit&\cr
& \!0.01\! &&\!0.000278548\!  && \!0.000278686\! && {}  &\cr
   height2pt&\omit&&\omit&&\omit&&\omit&&\omit&\cr
   \noalign{\hrule}
   height2pt&\omit&&\omit&&\omit&&\omit&&\omit&\cr
& 1 &&0.0234222  && 0.0245174 && 0.0245175  &\cr
   height2pt&\omit&&\omit&&\omit&&\omit&&\omit&\cr
   \noalign{\hrule}
   height2pt&\omit&&\omit&&\omit&&\omit&&\omit&\cr
& 100  &&0.13342 && 0.182541 && 0.182819  &\cr
     height2pt&\omit&&\omit&&\omit&&\omit&&\omit&\cr
 \noalign{\hrule}}}}}\hfill$$

The output $(z_k)$ of the algorithm in Section $\ref{s-03}$ is given in Table $10$.
\vskip-0.5truecm

 \begin{center}Table $10$.\quad  For different $b_4$, the
eigenvalue $\lz_0$, its lower bound $\dz_1^{-1}$ and $z_1, z_2$.
\end{center}
\vskip-0.4truecm

$$
\hfil\vbox{\hbox{\vbox{\offinterlineskip
 \halign{&\vrule#&\strut\quad\hfil#\hfil\quad&\vrule#&
 \quad\hfil#\hfil\quad&\vrule#&
 \quad\hfil#\hfil\quad\cr
    \noalign{\hrule}
   height2pt&\omit&&\omit&&\omit&&\omit&&\omit&\cr
   & {$\pmb{b_4}$}
   &&${\pmb{z_0=\dz_1^{-1}}}$ &&${\pmb{z_1}}$ &&${\pmb{z_2=\lz_0}}$ &\cr
   height2pt&\omit&&\omit&&\omit&&\omit&&\omit&\cr
   \noalign{\hrule}
   height2pt&\omit&&\omit&&\omit&&\omit&&\omit&\cr
& \!0.01\!  &&\!0.00027867\!  && \!0.000278686\! && {}  &\cr
   height2pt&\omit&&\omit&&\omit&&\omit&&\omit&\cr
   \noalign{\hrule}
   height2pt&\omit&&\omit&&\omit&&\omit&&\omit&\cr
& 1  &&0.0244003  && 0.024519 && 0.0245175  &\cr
  height2pt&\omit&&\omit&&\omit&&\omit&&\omit&\cr
   \noalign{\hrule}
   height2pt&\omit&&\omit&&\omit&&\omit&&\omit&\cr
& 100  &&0.179806 && 0.182912 && 0.182819  &\cr
     height2pt&\omit&&\omit&&\omit&&\omit&&\omit&\cr
   \noalign{\hrule}
& $10^6$  &&0.191917 && 0.195239 && 0.195145  &\cr
     height2pt&\omit&&\omit&&\omit&&\omit&&\omit&\cr
 \noalign{\hrule}
 }}}}\hfill$$
Once again, one sees the efficiency of our algorithm.}
\dexmp

Comparing the last two examples, especially their
$g$ and $\sqrt{x_i}$, it is obvious that
the latter is better than the former one. This suggests us
to choose the starting point $0$ carefully.
Here is an easier way to do so.
First, define a sequence $\{E_{\ell}\}$ of level sets as
follows. Let $E_0=\{N\}$ and $E_1=\{i\in E: a_{iN}>0\}$. At
the $k$th step, set
$$E_k=\{i\in E\setminus (E_0+\ldots +E_{k-1}): \exists j\in E_{k-1}
\text{ such that }a_{ij}>0\}.$$
The procedure should be stopped at $m$ if $E_{m+1}=\emptyset$.
Because of the irreducibility, each $i\in E$ should belong to
one of the level sets. Finally, regard one of $i_m\in E_m$
satisfying
$$a_{i_mj_{m-1}}=\min\{a_{ij}\!: i\in E_m,\; j\in E_{m-1}\}$$
as our initial $0$. However, for initial ${\tilde v}_0$, in practice, it is not necessary
to relabeling the states as we did in Example \ref{4-21}. What we
need is only replace the constraint $x_0=1$ by $x_{i_m}=1$
(at the same time, ``removing the first line'' is replaced by
``removing the $i_m$'s line'' in constructing the required
matrix) in solving $(x_i)$ without
change the original matrix $A$ or $Q$. One may need the relabeling
in computing $\dz_1$ defined in Section \ref{s-03}.

To conclude this subsection, we introduce a new construction of $z_0$ based on
$v_0$ defined by our general algorithm. It is an extension of $z_0=\dz_1^{-1}$
given in Section \ref{s-03}. To do so, we use $Q$, $(h_i)$ and $(x_i)$ defined
at the beginning of this subsection. Let ${\widetilde Q}_0$ be the matrix
obtained from
$${\widetilde Q}:=\text{\rm Diag}(h_i)^{-1}Q\,\text{\rm Diag}(h_i)$$
by modifying the last diagonal element ${\tilde q}_{N, N}$ so that the
sum of its last row becomes zero (i.e., removing the killing $c_N$). Next, let $\mu:=(\mu_0, \mu_1, \ldots, \mu_N)$
with $\mu_0=1$ be the solution to the equation
$$\mu {\widetilde Q}_0=0.$$
Since there are only $N$ variables $\mu_1, \ldots, \mu_N$, one may get the solution $\mu$
from the next equation
$${\widetilde Q}^{*\,\setminus \text{the last row}}\mu^*=0.$$
Here we remark that for a large class of $Q$-matrix $Q$, there is an explicit
representation of $\mu$ in terms of the non-diagonal elements of $Q$, refer to
\rf{cmf04}{Chapter 7}. Now, our new initial $z_0$ is defined to be $\dz_1^{-1}$:
\be \dz_1=\frac{1}{1-x_1}\max_{0\le n\le N}\bigg[\sqrt{x_n}\sum_{k=0}^n \mu_k \sqrt{x_k}+
\frac{1}{\sqrt{x_n}}\sum_{n+1\le j\le N}\mu_j x_j^{3/2}\bigg].\lb{a-11}\de

In contrast to the above examples which use only the automatic
$z_0=v_0^* A v_0$ (or $z_0=v_0^* (-Q) v_0$), here we use (\ref{a-11}).
Remember that this initial $z_0$ is for $-Q$, when we go back to the original
$A$, its initial becomes $\max_{i\in E}\sum_{j\in E} a_{ij} -z_0$.
The outputs of Examples \ref{x-18}--\ref{x-20} using $\dz_1^{-1}$ are listed
in Table $11$.
\vskip-0.5truecm

 \begin{center}Table $11$.\quad  The outputs of Examples \ref{x-18}--\ref{x-20} using $\dz_1^{-1}$.
\end{center}
\vskip-0.8truecm  

$$
\hfil\vbox{\hbox{\vbox{\offinterlineskip
 \halign{&\vrule#&\strut\quad\hfil#\hfil\quad&\vrule#&
 \quad\hfil#\hfil\quad&\vrule#&
 \quad\hfil#\hfil\quad\cr
    \noalign{\hrule}
   height2pt&\omit&&\omit&&\omit&&\omit&&\omit&\cr
   & {\pmb{Example}} && {${\pmb {z_0}}$}
   &&${\pmb{z_1}}$ &&${\pmb{z_2}}$ &&${\pmb{z_3}}$ &\cr
   height2pt&\omit&&\omit&&\omit&&\omit&&\omit&\cr
   \noalign{\hrule}
   height2pt&\omit&&\omit&&\omit&&\omit&&\omit&\cr
& \ref{x-18} && \!5.90016\! &&\!5.22268\!  && \!5.23611\! &&\!\!5.23607$=\lz_0$\!&\cr
   height2pt&\omit&&\omit&&\omit&&\omit&&\omit&\cr
   \noalign{\hrule}
   height2pt&\omit&&\omit&&\omit&&\omit&&\omit&\cr
& \ref{x-19} && 57.2719 &&36.236  && 36.2097 && 36.2094$=\lz_0$  &\cr
  height2pt&\omit&&\omit&&\omit&&\omit&&\omit&\cr
   \noalign{\hrule}
   height2pt&\omit&&\omit&&\omit&&\omit&&\omit&\cr
& \ref{x-20} && 30.4808 &&20.1294 && 24.8362 && 24.0778  &\cr
     height2pt&\omit&&\omit&&\omit&&\omit&&\omit&\cr
 \noalign{\hrule}}}}}\hfill$$
For Example \ref{x-20}, we need two more iterations: $z_4\approx 24.0294$, $z_5\approx 24.0293=\lz_0$.

Finally, we have an improved algorithm (for $Q$) as stated in Section \ref{s-03} (below Example \ref{t-09}) based on the use of $L^2(\mu)$ and the convex combination:
$$z_0=\xi \dz_1^{-1}+ (1-\xi) ( v_0, -Q v_0)_{\mu},\qqd \xi\in [0, 1].$$
The outputs of Examples \ref{x-18}--\ref{x-20} using the new $z_0$ with some $\xi$ are listed in Table $12$.

\begin{center}Table $12$.\quad  The outputs of Examples \ref{x-18}--\ref{x-20} using the new $z_0$ with $\xi\in (0, 1)$.
\end{center}
\vskip-0.8truecm  

$$
\hfil\vbox{\hbox{\vbox{\offinterlineskip
 \halign{&\vrule#&\strut\quad\hfil#\hfil\quad&\vrule#&
 \quad\hfil#\hfil\quad&\vrule#&
 \quad\hfil#\hfil\quad\cr
    \noalign{\hrule}
   height2pt&\omit&&\omit&&\omit&&\omit&&\omit&\cr
   &\!{\pmb{Example}}\!&&{$\pmb\xi$} && {${\pmb {z_0}}$}
   &&${\pmb{z_1}}$ &&${\pmb{z_2}}$ &&${\pmb{z_3=\lz_0}}$\! &\cr
   height2pt&\omit&&\omit&&\omit&&\omit&&\omit&\cr
   \noalign{\hrule}
   height2pt&\omit&&\omit&&\omit&&\omit&&\omit&\cr
& \ref{x-18}\!\! &&\!\!\!$1/3$\!&& \!5.04169\! &&\!5.24358\!  && \!5.23608\! &&\!\!{5.23607}\!&\cr
   height2pt&\omit&&\omit&&\omit&&\omit&&\omit&\cr
   \noalign{\hrule}
   height2pt&\omit&&\omit&&\omit&&\omit&&\omit&\cr
& \ref{x-19}&&$1/3$ && 35.4952 &&36.2657  && 36.2095 &&\!36.2094\!  &\cr
  height2pt&\omit&&\omit&&\omit&&\omit&&\omit&\cr
   \noalign{\hrule}
   height2pt&\omit&&\omit&&\omit&&\omit&&\omit&\cr
& \ref{x-20}&& $0.65$ && 24.0344 &&24.0161 && 24.0293 && {}  &\cr
     height2pt&\omit&&\omit&&\omit&&\omit&&\omit&\cr
 \noalign{\hrule}}}}}\hfill$$
From the table, it follows that the choice of $\xi$ depends on the model.
The safest case is $\xi=1$ and the dangerous one is $\xi=0$. When $\xi$
decreases in a neighborhood of $1$, the resulting $z_0$ can be usually improved.
Fortunately, for each class of models, such a $\xi$ can often be fixed.

This combination becomes more serious when $N$ is large since in that case
$( v_0, -Q v_0)_{\mu}$ is often an upper bound of $\lz_0$, which may be much
closer to other $\lz_j\ne\lz_0$ and so the algorithm would converge
to $\lz_j$ but not $\lz_0$.
Certainly, the convex combination idea is also meaningful for the first
two choices of $z_0$ introduced in the first subsection.

\section{Additional remarks and proofs}\lb{s-05}

In this section, we first prove a new result related to our earlier
study. Then we present some proofs of the results given in the last
two sections. Finally, we will make some remarks on the results
studied so far in the previous sections.

The next result solves an open question kept in our mind for
many years. For a given birth--death matrix $Q$ on $E$ with
$c_0=\ldots=c_{N-1}=0$ and $b_N:=c_N>0$, and a positive function
$f$ on $E$, define
$$I\!I(f)(i)=\frac{1}{f_i}\sum_{j=i}^N \frac{1}{\mu_j b_j}\sum_{k=0}^j\mu_k f_k,\qqd i\in E.$$

\prp\;\;{\cms For $Q$ and $I\!I$ given above, let $f_1 (>0\;\text{\cms on } E)$
be arbitrarily given function and define successively $f_{n+1}=f_n I\!I(f_n)$. Then this algorithm coincides with the inverse iteration given
in Lemma \ref{t-06} with $z=0$, even for infinite $N$. Furthermore, we have
$$\lz_0=\lz_{\min}(-Q)=\lim_{n\to\infty} I\!I(f_n)(i)^{-1}$$
for each $i\in E$. In particular, we have
$$\lim_{n\to\infty} \min_{i\in E} I\!I(f_n)(i)
=\frac{1}{\lz_0}=\lim_{n\to\infty} \max_{i\in E} I\!I(f_n)(i).$$}
\deprp

\prf
Consider the Poisson equation: $-Q f=g$ for a given $g$.
The solution is given by $f=gI\!I(g)$
(\rf{cmf10}{(2.7)--(2.9)}). It can be also written as
$f=(-Q)^{-1}g$. By setting $g=f_1$ and $f=f_2$, it follows that
$$f_2=(-Q)^{-1} f_1= f_1 I\!I(f_1).$$
Now, by iteration, we get
$$f_{n+1}=(-Q)^{-n} f_1= f_n I\!I(f_n),\qqd n\ge 1.$$
We have thus proved the first assertion. Therefore,
$$I\!I(f_n)=\frac{f_{n+1}}{f_n}=\frac{(-Q)^{-n} (f_1)}{(-Q)^{-n+1} (f_1)}
\to \frac{1}{\lz_0}\qqd \text{as }n\to \infty$$
by the last assertion of Lemma \ref{t-06} with $z=0$. The last assertion of
the proposition then follows since on a finite set, the pointwise convergence
implies the uniform one.
\deprf

We remark that the last proposition is meaningful once the Poisson
equation $-Q f =g$ is solvable. In parallel, Lemma \ref{t-06} improves
the approximating procedures studied in \ct{cmf10} and related
publications.

\medskip

Now, we turn to prove Proposition \ref{t-08-1} and Corollary \ref{t-11-1}.

\nnd{\bf Proof of Proposition \ref{t-08-1}}.\;\; (a) First, we follow the setup and notation in \ct{chzh} (where a more general situation is studied) for a moment. Define
$$\begin{gathered}
M_{N-1}(h)= \sum_{k=0}^{N-1}{\tilde q}_N^{(k)}\sum_{j=0}^{k}\frac{{\widetilde F}_k^{(j)}h_j}{q_{j,\, j+1}},\\
N_n(h)= \sum_{k=0}^{n}\sum_{j=0}^{k}\frac{{\widetilde F}_k^{(j)}h_j}{q_{j,\, j+1}},
\qqd 0\le n< N.\end{gathered}$$
Then the solution given in \rf{chzh}{Proposition 2.6} can be rewritten as
$$g_n =\frac{f_N+M_{N-1}(f)}{c_N+M_{N-1}(c_{\cdot})}
\big[1-N_{n-1}(c_{\cdot})\big]+N_{n-1} (f),\qd N_{-1}:=0,\qqd 0\le n\le N.$$
By an exchange of the order of the summations, we can rewrite $M_n$
and $N_n$ as follows:
$$\begin{gathered}
M_{N-1}(h)= \sum_{j=0}^{N-1}\frac{h_j}{q_{j,\,j+1}}\sum_{k=j}^{N-1}{\tilde q}_N^{(k)}{\widetilde F}_k^{(j)},\\
N_n(h)= \sum_{j=0}^{n} \frac{h_j}{q_{j,\, j+1}}\sum_{k=j}^{n}{\widetilde F}_k^{(j)},
\qqd 0\le n< N.\end{gathered}$$
Here for finite $N$, the element $q_{N, N+1}$ is replaced by $c_N$ by our convention. Thus, by \rf{chzh}{(1.1)}, we get
$$M_{N-1}(h)= c_N\sum_{j=0}^{N-1}\frac{h_j}{q_{j,\,j+1}}{\widetilde F}_N^{(j)}.$$
Since by \rf{cmf14}{Proposition 4.1}, we have ${\widetilde F}_{i+m}^{(i)}=G_{m, m}^{(i)}$. It follows that
$$\begin{gathered}
M_{N-1}(h)= c_N\sum_{j=0}^{N-1}\frac{h_j}{q_{j,\,j+1}} G_{N-j,\, N-j}^{(j)},\\
N_n(h)= \sum_{j=0}^{n} \frac{h_j}{q_{j,\, j+1}}\sum_{k=0}^{n-j} G_{k,\,k}^{(j)},
\qqd 0\le n< N.\end{gathered}$$
Applying this solution to the birth--death context and setting
$f=-v$, $g=w$, replacing the original $c_{\cdot}$ used in \ct{chzh} by
$z-c_{\cdot}$, we obtain
$$g_n =\frac{-v_N-M_{N-1}(v)}{z-c_N+M_{N-1}(z-c_{\cdot})}
\big[1-N_{n-1}(z-c_{\cdot})\big]-N_{n-1} (v),\qqd 0\le n\le N.$$
Equivalently,
$$g_n =\frac{v_N+M_{N-1}(v)}{c_N-z+M_{N-1}(c_{\cdot}-z)}
\big[1+N_{n-1}(c_{\cdot}-z)\big]-N_{n-1} (v),\qqd 0\le n\le N.$$
This gives us the required conclusion.
\deprf

\medskip

\nnd{\bf Proof of Corollary \ref{t-11-1}}.\;\;The proof is quite straightforward. Choose $m$ large enough such that
$$A:=m I +Q$$
is a nonnegative matrix. Then  $-Q=m I -A$. Hence
$\lz_0(Q)=m -\rho(A)$. The proof now is a direct application of
the Collatz--Wielandt formula:
$$\aligned
m -\rho(A)&= m -\inf_{x>0}\max_i \frac{(Ax)_i}{x_i}
  =  \sup_{x>0}\min_i \frac{(-Qx)_i}{x_i},\\
m -\rho(A)&= m -\sup_{x>0}\min_i \frac{(Ax)_i}{x_i}
=  \inf_{x>0}\max_i \frac{(-Qx)_i}{x_i}.\qquad \square
\endaligned$$

\medskip

It is now ready to make some additional remarks on the results in the
previous sections. The two algorithms as well as their convergence and the Collatz--Wielandt formula can
be found easily from Wikipedia. From which, one knows that
the Power Iteration was first appeared in 1929 \ct{mpg29} and the Inverse
Iteration appeared in 1944 \ct{hw44}. These algorithms are taught
for undergraduate students on the course of computations and are included in many books, see for instance
\ct{mc00, jhw65, gl13}. In particular, the Appendix at the end of
Section \ref{s-03} is modified from \rf{gl13}{pages 584--585}.

We now say a few words about the unusual word ``complete'' used at the
end of the first section for the results obtained in Section \ref{s-03}.
Actually, this is one of the 16 situations with $N\le\infty$ we have worked out so far to
have a unified estimation of the principal eigenvalue:
\be (4\dz)^{-1}\le \dz_1^{-1}\le \lz_0\le {\dz_1'}^{-1}\le \dz^{-1}\lb{5-00} \de
for some constants $\dz, \dz_1, \dz_1'$, where $\dz_1$ is the one we
have used in Section \ref{s-03} for the initial $z_0$. Besides, we often have in practice that $1\le \dz_1/\dz_1'\le 2$. Thus, the efficiency
of the initial $(v_0, z_0)$ introduced in Section \ref{s-03} comes with
no surprising. More precisely, the initial $(v_0, z_0)$ is taken from
the first step of our approximating procedure:
\rf{cmf10}{Theorem 3.3\,(1) and (3.4)}. Example \ref{t-01} here is a truncated one from
\rf{cmf10}{Example 3.6} where $N=\infty$, $\lz_0=1/4$ and $\dz_1=4$ which is sharp.
Certainly, this is still not enough to claim that we can arrive at such
a precise approximation in the second iteration.
The story on the estimation of the principal eigenvalue, or more
general on the estimation of the stability speed is too long to talk here
and so the author is planning to publish a survey article \ct{cmf16}.
For earlier progress, refer to \ct{cmf05} which includes a lot of information up to 2004,
or a more recent paper \ct{cmf10}.

Next, we discuss the sequence $(h_0, h_1, \ldots, h_N)$ used in Sections
\ref{s-03} and \ref{s-04}. The role of the sequence is to keep the same
spectrum of the original $Q$ and its $H$-transform
${\widetilde Q}$:
\be {\widetilde Q}=\text{\rm Diag}(h_i)^{-1}Q\,\text{\rm Diag}(h_i).\lb{4-8-2}\de
Certainly, $Q$ and ${\widetilde Q}$ have the same diagonals. Next, define
\be P=(p_{ij}: i,j\in E):=\text{\rm Diag}(q_i^{-1})\,{\widetilde Q}+I,\lb{4-8-3} \de
which is the matrix used in Section \ref{s-04}.
Note that even though the sequence $(c_i)$ in the original $Q$ can be non-zero,
the resulting ${\tilde c}_k=0$ for every $k<N$ but ${\tilde c}_N>0$
for the matrix ${\widetilde Q}$.
For a given measure $\mu$, set ${\tilde \mu}=h^2\mu$ (i.e., ${\tilde \mu}_i=h_i^2\mu_i$ for each $i\in E$), the transform ${\tilde f}=f/h$
gives us an isometry between $L^2(\mu)$ and $L^2({\tilde\mu})$ and then
an isospectrum of $Q$ on $L^2(\mu)$ and ${\widetilde Q}$ on $L^2({\tilde\mu})$.
This technique is due to \ct{chzhx}. See also \ct{cmf14}.
Now, if $\tilde g$ is an approximating eigenvector corresponding to $\tilde\lz_0$ of $\widetilde Q$, then, $g:=h {\tilde g}$ is an
approximating eigenvector corresponding to $\lz_0$ of $Q$, due to the isospectral
property of $Q$ and $\widetilde Q$.
Because
$$\|\tilde g\|_{L^2(\tilde\mu)}=\|g\|_{L^2(\mu)},\qqd
\big({\tilde g}, {\widetilde Q}{\tilde g}\big)_{\tilde\mu}
=({g}, {Q}{g})_{\mu},$$
by \ct{chzhx}, we have
\be \frac{\big({\tilde g}, -{\widetilde Q}{\tilde g}\big)_{\tilde\mu}}{\|\tilde g\|_{L^2(\tilde\mu)}}=
\frac{({g}, {-Q}{g})_{\mu}}{\|g\|_{L^2(\mu)}}
=\frac{g^* (-Q) g}{\sqrt{g^*g}}, \lb{5-13}  \de
here we assume that $\mu_k\equiv 1$ for simplicity. This means that we can
estimate the maximal eigenpair $(\lz_0, g)$ of $Q$ in terms of
the one $\big(\tilde\lz_0, \tilde g\big)$ of ${\widetilde Q}$. More precisely, the maximal eigenvalue ${\tilde g}$ of
$\widetilde Q$ is approximated by $\fz$ in the context of
Section \ref{s-03} (or by $x=(x_i)$ in Section \ref{s-04}). Now, in Section \ref{s-03} for instance, ${\tilde v}_0=h\sqrt{\fz}$ is an approximation of
the maximal eigenvector $g$ of $Q$. With $v_0={\tilde v}_0/\sqrt{{\tilde v}_0^* v_0}$, equation (\ref{5-13}) leads to our first approximation of $\lz_0$:
$$v_0^* (-Q) v_0=z_0.$$

Now, our task is to show that the sequence $(x_i)$ defined in Section \ref{s-04} is an extension of $(\fz_i)$ given in Section \ref{s-03}. To this end,
recall that the matrix ${\widetilde Q}$ defined by (\ref{4-8-2}) is again a $Q$-matrix.
Hence, the matrix $P=(p_{ij}: i,j\in E)$ defined by (\ref{4-8-3})
is just the embedding chain of ${\widetilde Q}$.
Note that here $p_{ii}=0$ for each $i\in E$. By the construction of $(h_i)$,
we have $\sum_{j\in E}p_{ij}=1$ for each $i\le N-1$ but $\sum_{j\in E}p_{Nj}^{}<1$,
refer to \ct{chzhx}. The equation for $(x_i)$ in (\ref{4-8-1})
 can be rewritten as
\be x_n=\sum_{j\in E}p_{ij} x_j,\qqd 1\le n\le N,\; x_0=1.\lb{5-11}\de
In probabilistic language, the solution $(x_i)$ (or the minimal
solution $(x_i^*)$ when $N=\infty$) to equation (\ref{5-11}) is the
probability of first hitting $0$ of the $Q$-process with $Q$-matrix
$\widetilde Q$ or its embedding sub-Markov chain with transition matrix $P=(p_{ij})$, starting from $i$. Refer to \rf{cmf04}{Lemma 4.46}.

We are now going to prove the following result.

\lmm\lb{t-21}\;\;{\cms For birth--death matrix, the solution $(x_i)$ to equation $(\ref{5-11})$
coincides with $(\fz_i)$ (up to a constant) used in Section $\ref{s-03}$.}
\delmm

Before prove Lemma \ref{t-21}, let us discuss the relation of these sequence
with the recurrence of the Markov chain in the case of $N=\infty$. First, it is known by
\rf{cmf04}{Theorem 4.55\,(1) and the second line of page 161} that
a birth--death process is recurrent iff
$$b_0\sum_{n=1}^{\infty} \frac{a_1\cdots a_n}{b_1\cdots b_n}=b_0\sum_{n=1}^{\infty} \frac{1}{\mu_n b_n}
=\infty.$$
For simplicity, set
$$F_n^{(0)}=\frac{a_1\cdots a_n}{b_1\cdots b_n},\qqd n\ge 1.$$
The sequence $\big\{F_n^{(0)}\big\}_{n\ge 1}$ is a very special case of $\big\{{\widetilde F}_n^{(j)}\big\}_{n\ge 1}$ used in the proof of Proposition \ref{t-08-1}. Refer to
\ct{chzh} and \rf{cmf04}{\S 4.5} for more details.
 Note that $(\fz_n)$ is just the tail series of $\sum_{n=1}^{\infty} F_n^{(0)}$ provided $N=\infty$. On the other hand, by
\rf{cmf04}{Lemma 4.46}, the process is recurrent iff the minimal solution
$(x_i^*)$ to the equation (\ref{5-11}):
$$x_n=\frac{b_n}{a_n+b_n}x_{n+1}+\frac{a_n}{a_n+b_n}x_{n-1},\qqd n\ge 1,\; x_0:=1$$
is equal to one identically. Rewrite the equation as
$$x_n-x_{n+1}= \frac{a_n}{b_n}(x_{n-1}-x_n),\qqd n\ge 1.$$
By induction, it follows that
$$x_n-x_{n+1}=F_n^{(0)}(x_0-x_1),\qqd n\ge 1.$$
Hence
$$\aligned
x_n-x_{N+1}&=(x_0-x_1)\sum_{k=n}^N F_k^{(0)},\qd
x_1-x_n=(x_0-x_1)\sum_{k=1}^{n-1} F_k^{(0)},\qd n\ge 1.
\endaligned
$$
Equivalently,
$$\aligned
x_n-x_{N+1}&=(x_0-x_1)\sum_{k=n}^N F_k^{(0)},\qd
x_0-x_n=(x_0-x_1)\sum_{k=0}^{n-1} F_k^{(0)},\qd n\ge 0.
\endaligned
$$
since $F_0^{(0)}=1$ by convention. If $\sum_{k=0}^{\infty} F_k^{(0)}=\infty$, then from
the second equation, we must have $x_1\!=\!1$ (since $x_0\!=\!1$) and then have the
unique solution $x_i\equiv 1$. Therefore, the minimal solution
$x_i^*\equiv 1$ and so the process is recurrent. Conversely, if
$\sum_{k=0}^{\infty} F_k^{(0)}<\infty$, then from the first equation
above, we obtain
$$x_0-x_1= \frac{x_0-x_{\infty}}{\sum_{j=0}^{\infty} F_j^{(0)}},$$
and then
$$x_n-x_{\infty}= \frac{x_0-x_{\infty}}{\sum_{j=0}^{\infty} F_j^{(0)}}
\sum_{k=n}^{\infty} F_k^{(0)},\qqd n\ge 0.$$
Equivalently,
$$x_n= \frac{\sum_{k=n}^{\infty} F_k^{(0)}}{\sum_{j=0}^{\infty} F_j^{(0)}}
+ x_{\infty}\frac{\sum_{k=0}^{n-1} F_k^{(0)}}{\sum_{j=0}^{\infty} F_j^{(0)}},
\qqd n\ge 0.$$
Clearly, for each given $x_{\infty}\in [0, 1]$, using this formula,
we obtain a solution $(x_n)$ to the equation. Thus, the minimal solution
should be as follows:
$$x_n^*= \frac{\sum_{k=n}^{\infty} F_k^{(0)}}{\sum_{j=0}^{\infty} F_j^{(0)}},
\qqd n\ge 0,$$
which is clearly less than one for $n\ge 1$ since $\sum_{j=0}^{\infty} F_j^{(0)}<\infty$.
\medskip

\nnd{\bf Proof of Lemma \ref{t-21}}\;\; For finite state $\{0, 1, \ldots, N\}$, since there is a killing $b_N>0$,
the minimal solution is as follows:
$$x_n^*= \frac{\sum_{k=n}^{N} F_k^{(0)}}{\sum_{j=0}^{N} F_j^{(0)}},
\qqd n= 0, 1, \ldots, N.$$
In other words, up to a constant, we have
$$ \fz_n= \sum_{k=n}^{N} F_k^{(0)}= \frac{1}{1-x_1^*} x_n^*, \qqd n= 0, 1, \ldots, N.$$
That is what we required. \deprf

Finally, we remark that the story for one-dimensional diffusions should be in parallel to Section \ref{s-03}. The algorithm presented in Section
\ref{s-04} may not be complete since the lack of an analog of (\ref{5-00}).

\subsection*{Summary} \quad This paper deals with the efficient initials
for the Rayleigh quotient iteration. Here are suggestions for the use of the results in the previous sections of the paper on computing the maximal eigenpair.

(1) If the iterations are easy (small size of $A$, for instance), one simply
adopts the simplest algorithm: Section \ref{s-4-1} with Choice I, or more
effectively with the convex combination of Choice I and Choices II:
$$z_0=\xi \max_{i\in E}A_i+(1-\xi) v_0^* A v_0, \qqd \xi\in [0, 1].$$   
More especially, $\xi=7/8$ for instance. Certainly, one may use the choice III for $z_0$.

(2) If the given matrix is nearly tridiagonal (after a suitable relabeling if necessary) or the Lanczos tridiagonalization procedure is
suitable, one use the method introduced in
Section \ref{s-03}. The computation there is rather explicit and it works even
for $N=\infty$.

(3) In general, one uses the algorithm given in Section \ref{s-4-2}. Note that
at each step of the Rayleigh quotient iteration, one has to solve a linear
equation. Here for the initials, we have to solve two more linear
equations.

\section{Next to maximal eigenpair}\lb{s-06}

After an earlier version of the paper containing the first five sections was submitted,
the author found a natural way to study the next to the maximal eigenpair. In this section, we restrict
ourselves to the easier
case that $A_i:=\sum_{j\in E} a_{ij}$ is a constant. Then the maximal
eigenpair is simply $(A_0, \bbb{1})$ (where $\bbb{1}$ is the constant function
having value $1$ everywhere), as mentioned before. By a shift
if necessary, we return to the problem for a $Q$-matrix which is especially
valuable since its next eigenvalue describes the ergodic rate of the
corresponding Markov chain. In this setup, the minimal eigenpair
$(\lz_0=0,\, g_0=\bbb{1})$ of $-Q$ is known and we are looking for the next eigenpair
$(\lz_1, g_1)$. Clearly, $g_1$ should be orthogonal to $g_0$ in $L^2(\pi)$-sense
for the stationary distribution $\pi$ of the process corresponding to
the given matrix $Q$. This is the reason why we often use $v-\pi v$ in what follows
for constructing a mimic of the eigenvector $g_1$. Besides, we need the assumption
that $\lz_1>|\lz_j|$ for every $j> 1$ to guarantee the convergence of our algorithms.

Once again, let us begin our study with
a tridiagonal conservative $Q$-matrix
$$Q\!=\!\left(\!\begin{array}{ccccc}
-b_0\!\!\! & b_0 &0&0 &\cdots \\
a_1 &\!\!\! -(a_1 + b_1) & b_1 &0 &\cdots \\
0& a_2 &\!\!\! -(a_2 + b_2) & b_2 &\cdots \\
\vdots &\vdots &\ddots &\ddots &\ddots \\
0& 0 &\qquad 0 &\quad a_N^{}\;\; & -a_N^{}
\end{array}\!\right)\!,$$
where $a_i> 0,\; b_i> 0$.
Define $(\mu_k: k\in E)$ as in \S \ref{s-03}. Then we have the probability
distribution $\pi=(\pi_0, \pi_1,\ldots, \pi_N)$: $\pi_k=\mu_k/\sum_{j\in E}\mu_j$.
Again, denote by $(\cdot, \cdot)_{\mu}$ and $\|\cdot\|_{\mu}$ the inner product
and norm in $L^2(\mu)$, respectively. Next, set
$$\fz_n=\sum_{j\le n-1}\frac{1}{\mu_j b_j},\qqd n\in E.$$
To define our initial $v_0$, let
$${\tilde v}_0=\big(\sqrt{\fz_0},\, \sqrt{\fz_1},\, \ldots,\, \sqrt{\fz_N}\,\big)^*,
\qqd {\bar v}_0={\tilde v}_0-\pi {\tilde v}_0.$$
We can now introduce our algorithm in the present situation as follows. Choose initials
\be v_0=\frac{{\bar v}_0}{\|{\bar v}_0\|_{\mu}},\qqd
z_0=\frac{({\bar v}_0, -Q{\tilde v}_0)_{\mu}}{\|{\bar v}_0\|_{\mu}^2}.\lb{6-17} \de
At the $k\,(k\ge 1)$th step, let $w_{k}$ be the solution to the equation
$$(-Q-z_{k-1})w_{k}=v_{k-1}$$
and set
$$v_k=\frac{w_k}{\|w_k\|_{\mu}},\qqd z_k=(v_k, -Q v_k)_{\mu}.$$

   We remark that here in defining $v_k\,(k\ge 1)$, we do not need
to use $w_k - \pi w_k.$ The reason is as follows. If $\pi  v=0$ and $w$ solves the equation
   $$(-Q-z)w=v$$
for some constant $z\ne 0$,   then
$$0=\pi v=\pi (-Q-z)w=-z \pi  w, $$
and so $\pi  w=0$. Therefore, we have $\pi  w_k=0$ for each $k\ge 1$
since so does the initial $v_0$: $\pi v_0=0$.

Instead of $z_0$ given in (\ref{6-17}), there is another choice. Define
$$\ez_1=\max_{0\le i\le N-1}\frac{1}{\mu_i b_i \big[{\tilde v}_0(i+1)
- {\tilde v}_0(i)\big]}\sum_{j=i+1}^N \mu_j {\bar v}_0(j).$$
Then one may choose
\be z_0=\ez_1^{-1}\lb{6-18} \de
as an initial.

Here the initials $\tilde v_0$ and $z_0$ are taken from \rf{cmf01}{Theorem 2.2\,(1)} or \rf{cmf05}{Theorem 1.5\,(2)}.
Certainly, we can adopt the convex combination of those given in
(\ref{6-17}) and (\ref{6-18}):
\be z_0= \xi \ez_1^{-1}+ (1-\xi) {({\bar v}_0, -Q{\tilde v}_0)_{\mu}}{\|{\bar v}_0\|_{\mu}^{-2}},\qqd \xi \in [0, 1].\lb{6-18-1}\de

We now consider an example modified from Example \ref{t-01}.
\xmp\lb{t-30} \;\;{\cms Let $E=\{0, 1, \ldots, 7\}$ and
$$Q=\left(
\begin{array}{cccccccc}
 -1 & 1 & 0 & 0 & 0 & 0 & 0 & 0 \\
 1 & -5 & 2^2 & 0 & 0 & 0 & 0 & 0 \\
 0 & 2^2 & -13 & 3^2 & 0 & 0 & 0 & 0 \\
 0 & 0 & 3^2 & -25 & 4^2 & 0 & 0 & 0 \\
 0 & 0 & 0 & 4^2 & -41 & 5^2 & 0 & 0 \\
 0 & 0 & 0 & 0 & 5^2 & -61 & 6^2 & 0 \\
 0 & 0 & 0 & 0 & 0 & 6^2 & -85 & 7^2 \\
 0 & 0 & 0 & 0 & 0 & 0 & 7^2 & -7^2
\end{array}
\right).$$
Then we have $\mu_k\equiv 1$, $\lz_1(Q)\approx 0.820539$ with eigenvector
$$\approx(-3.95053, -0.708966, 0.246859, 0.649164, 0.842169, 0.93805, 0.983254, 1)^*.$$
Starting from ${\bar v}_0$:
$$(-4.79299,\, -0.0815238,\, 0.474589,\, 0.70372,\, 0.828504,\, 0.906932,\, 0.960767,\, 1)^*,$$
for different initial $z_0$, the outputs are given in Table 13.

\begin{center}Table 13.\quad Outputs for different initial $z_0$ (Example \ref{t-30}).\end{center}
\vskip-0.5truecm

$$
\hfil\vbox{\hbox{\vbox{\offinterlineskip
 \halign{&\vrule#&\strut\quad\hfil#\hfil\quad&\vrule#&
 \quad\hfil#\hfil\quad&\vrule#&
 \quad\hfil#\hfil\quad\cr
    \noalign{\hrule}
   height2pt&\omit&&\omit&&\omit&&\omit&\cr
   & \pmb{Choice} && $\pmb {z_0}$ && $\pmb {z_1}$ &&  $\pmb{z_2=\lz_1}$ &\cr
   height2pt&\omit&&\omit&&\omit&&\omit&\cr
   \noalign{\hrule}
   height2pt&\omit&&\omit&&\omit&&\omit&\cr
   & (\ref{6-17}) && $0.902633$ && $ 0.820614$ && $0.820539$
  &\cr
   height2pt&\omit&&\omit&&\omit&&\omit&\cr
   \noalign{\hrule}
   height2pt&\omit&&\omit&&\omit&&\omit&\cr
& (\ref{6-18}) && $0.456343$ && $0.8216$ && $0.820539$
    &\cr
   height2pt&\omit&&\omit&&\omit&&\omit&\cr
   \noalign{\hrule}
   height2pt&\omit&&\omit&&\omit&&\omit&\cr
   & (\ref{6-18-1}) && $0.724117$ && $0.820629$ && $0.820539$
    &\cr
   height2pt&\omit&&\omit&&\omit&&\omit&\cr
 \noalign{\hrule}}}}}\hfill$$
We remark that for this and the next example, the parameter
$\xi$ in (\ref{6-18-1}) is specified to be $2/5$.
}\dexmp

The next example has non-trivial $(\mu_k)$.       
\xmp\lb{t-31}{\cms Let
$$Q=\left(
\begin{array}{ccccc}
-5 & 5 & 0 & 0 & 0\\
3 & -7 & 4 & 0 & 0\\
0 &  2 & -3 & 1 & 0\\
0 & 0 & 10 & -16 & 6\\
0 & 0 & 0 & 11 & -11
\end{array}
\right).$$
Then
$$\mu_0=1,\;\; \mu_1= 5/3,\;\;
\mu_2= 10/3,\;\; \mu_3= 1/3,\;\; \mu_4= 2/11.$$
The eigenvalues of $-Q$ are as follows:
$$22.348,\;\; 10.6857,\;\; 5.92951,\;\; 3.03673,\;\; 0.$$
With
$${\tilde v}_0=\frac 1 {2\sqrt{5}} \big(0,\; 2,\; \sqrt{7}, \sqrt{13},\; \sqrt{23}\,\big),$$
for different initial $z_0$, the outputs are given in Table 14.

\begin{center}Table 14.\quad Outputs for different initial $z_0$ (Example \ref{t-31}).\end{center}
\vskip-0.5truecm

$$
\hfil\vbox{\hbox{\vbox{\offinterlineskip
 \halign{&\vrule#&\strut\quad\hfil#\hfil\quad&\vrule#&
 \quad\hfil#\hfil\quad&\vrule#&
 \quad\hfil#\hfil\quad\cr
    \noalign{\hrule}
   height2pt&\omit&&\omit&&\omit&&\omit&\cr
   & \pmb{Choice} && $\pmb {z_0}$ && $\pmb {z_1}$ &&  $\pmb {z_2=\lz_1}$ &\cr
   height2pt&\omit&&\omit&&\omit&&\omit&\cr
   \noalign{\hrule}
   height2pt&\omit&&\omit&&\omit&&\omit&\cr
   & (\ref{6-17}) && $3.84977$ && $ 3.05196$ && $3.03673$
  &\cr
   height2pt&\omit&&\omit&&\omit&&\omit&\cr
   \noalign{\hrule}
   height2pt&\omit&&\omit&&\omit&&\omit&\cr
& (\ref{6-18}) && $1.72924$ && $3.05715$ && $3.03673$
    &\cr
   height2pt&\omit&&\omit&&\omit&&\omit&\cr
   \noalign{\hrule}
   height2pt&\omit&&\omit&&\omit&&\omit&\cr
   & (\ref{6-18-1}) && $3.00156$ && $3.03675$ && $3.03673$
    &\cr
   height2pt&\omit&&\omit&&\omit&&\omit&\cr
 \noalign{\hrule}}}}}\hfill$$
}\dexmp

Next, consider the general conservative $Q$-matrices $Q=(q_{ij}: i,j\in E)$. Here the
conservativity means that $\sum_{j\in E} q_{ij}=0$ for every $i\in E$. Next, define
an auxiliary $Q$-matrix $Q_1$ which coincides with $Q$ except replacing the element $q_{NN}$ by $c\, q_{NN}$, where $c>1$ is an arbitrary constant and is fixed to be 1000 in what follows for simplicity.

Following \S 4 (replacing $Q$ by $Q_1$), let $(x_0, x_1, \ldots, x_N)$ (with $x_0=1$)
be the solution to the equation
\be x^{\setminus \text{$0$'s row}}
=P^{\setminus \text{$0$'s row}}\,x,    \de
where
$$P=\text{\rm Diag}\big(q_0^{-1},\, q_1^{-1}, \,\ldots,\, q_{N-1,N-1}^{-1},\,(c q_{NN})^{-1}\big)Q_1+I.$$
To go further, we need $\mu=(\mu_0, \mu_1, \ldots, \mu_N)$ with $\mu_0=1$, which is the
same as defined in \S \ref{s-04}: the solution to the equation
$${Q}^{*\,\setminus \text{the last row}}\mu^*=0.$$
Having $x$ and $\mu$ at hand, we are ready to define our initials. For each $r\in [0, 1]$, to be specified later, define
\begin{align}
&{\tilde v}_0=\big(r, \,\sqrt{1-x_1},\,\ldots, \,\sqrt{1-x_N}\,)^*;\qqd {\bar v}_0={\tilde v}_0-\mu{\tilde v}_0\bigg/\sum_{k=0}^N \mu_k,  \nonumber\\
&v_0=\frac{{\bar v}_0}{\|{\bar v}_0\|_{\mu}},\qqd
z_0=\frac{({\bar v}_0, -Q {\tilde v}_0)_{\mu}}{\|{\bar v}_0\|_{\mu}^2}.  \lb{6-20}
\end{align}
Because ${\tilde v}_0$ depends on $r$, so do ${\bar v}_0$, $v_0$ and $z_0=: z_0(r)$.
Choose $r_0\in [0, 1]$ so that
$$z_0(r_0)\approx \inf_{r\in [0, 1]} z_0(r).$$
Corresponding to this specified $r_0$, we obtain our initials $v_0$ and $z_0$.
This minimizing procedure in $r$ is necessary for avoiding collapse since we are
in a more sensitive situation than before. Then
the iteration procedure is exactly the same as we used several times before.

The reason we adopt a large $c=1000$ here is that for a larger $c$, its minimal
eigenvalue $\lz_0(Q_1)$ is closer to, but less than, the eigenvalue $\lz_1(Q)$ we are interested. Refer to  \rf{cmf00}{Proposition 3.2} for more details. Thus, one may regard the former as an approximation of the latter.
In other words, we can use an alternative initial
\be z_0=\lz_0(Q_1)\;\text{or its estimates studied in previous sections.} \lb{6-21}\de

Certainly, one can define a convex combination of those given in (\ref{6-20})
and (\ref{6-21}) in an obvious way, but it is omitted here. The use of $\lz_0(Q_1)$ seems necessary (especially for large $N$) to avoid some pitfall, as mentioned before.

The next example is interesting for which some of its eigenvalues are complex but the one we are interested is real.
\xmp\;{\cms Let
$$Q =\begin{pmatrix}
-30 & 30 &  0 &  0\\
1/5 &  -17 & 84/5 &  0\\
11/28 & 275/42 & -20 & 1097/84\\
55/3291 & 330/1097 & 588/1097 & -2809/3291\end{pmatrix}.$$
Then
$$Q_1 =\begin{pmatrix}
-30 & 30 &  0 &  0\\
1/5 &  -17 & 84/5 &  0\\
11/28 & 275/42 & -20 & 1097/84\\
55/3291 & 330/1097 & 588/1097 & -2809000/3291\end{pmatrix}.$$
The eigenvalues of $-Q$ and $-Q_1$ are
$$
29.8411 + 2.45214\,i,\;\; 29.8411 - 2.45214\,i,\;\; 8.17131,\;\; 0
$$
and
$$853.548, \;\; 29.8249 + 2.46241\,i,\;\; 29.8249 - 2.46241\, i,\;\; 7.34195$$
respectively. Using (\ref{6-20}) with $r_0\approx 0.951$, the output is
$$z_0\approx 7.73667,\;\; z_1\approx 8.15021,\;\; z_2\approx 8.17129,\;\; z_3\approx 8.17131.$$
While using (\ref{6-21}), the output is
$$z_0\approx 7.34195,\;\; z_1\approx 8.13216,\;\; z_2\approx 8.17124,\;\; z_3\approx 8.17131.$$
}\dexmp

Here is one more example.
\xmp\;{\cms Let
$$Q =\begin{pmatrix}
-57 & 118/27 & 91/9 & 1148/27\\
135/59 & -52 & 637/59 & 2296/59\\
243/91 & 590/91 & -47 & 492/13\\
351/287 & 118/41 & 195/41 & -62/7\end{pmatrix}.$$
Then
$$Q_1 =\begin{pmatrix}
-57 & 118/27 & 91/9 & 1148/27\\
135/59 & -52 & 637/59 & 2296/59\\
243/91 & 590/91 & -47 & 492/13\\
351/287 & 118/41 & 195/41 & -62000/7\end{pmatrix}.$$
The eigenvalues of $-Q$ and $-Q_1$ are
$$
59.3118,\;\; 58,\;\; 47.5454,\;\; 0
$$
and
$$8857.18, \;\; 59.2467,\;\; 58,\;\; 38.7143$$
respectively. Using (\ref{6-20}) with $r_0\approx 0.953$, the output is
$$z_0\approx 47.5318,\;\; z_1\approx 47.5453,\;\; z_2\approx 47.5454.$$
While using (\ref{6-21}), the output is
$$z_0\approx 38.7143,\;\; z_1\approx 47.5343,\;\; z_2\approx 47.5453,\;\; z_3\approx 47.5454.$$
}\dexmp

\medskip

\nnd{\bf Acknowledgments}. {\small
The main results of the paper have been reported
at Anhui Normal University,
Jiangsu Normal University, the International Workshop on SDEs and Numerical Methods
at Shanghai Normal University, Workshop on Markov Processes and Their Applications at
Hunan University of Arts and Science, and Workshop of Probability Theory with Applications at University of Macau.
The author acknowledges Professors Dong-Jin Zhu, Wan-Ding Ding,
Ying-Chao Xie, Xue-Rong Mao, Xiang-Qun Yang, Xu-Yan, Xiang, Jie Xiong, Li-Hu Xu, and their teams for very warm hospitality and financial
support. The author also thanks Ms Yue-Shuang Li for her assistance in computing large matrices.
Research supported in part by the
         National Natural Science Foundation of China (No. 11131003),
         the ``985'' project from the Ministry of Education in China,
and the Project Funded by the Priority Academic Program Development of
Jiangsu Higher Education Institutions.
}

\medskip

\nnd {\small
Mu-Fa Chen\\
School of Mathematical Sciences, Beijing Normal University,
Laboratory of Mathematics and Complex Systems (Beijing Normal University),
Ministry of Education, Beijing 100875,
    The People's Republic of China.\newline E-mail: mfchen@bnu.edu.cn\newline Home page:
    http://math.bnu.edu.cn/\~{}chenmf/main$\_$eng.htm}



\begin{thebibliography}{10}
\setlength{\itemsep}{-0.8ex}
{\small

\bibitem{cmf00}
Chen, M.F. (2000)
{\it  Explicit bounds of the first eigenvalue}.
Sci. China {\rm (A)} 43(10), 1051--1059.\lb{cmf00}

\bibitem{cmf01}
Chen, M.F. (2001)
{\it Variational formulas and approximation theorems for the first
  eigenvalue}.
Sci. China {\rm (A)} 44(4), 409--418. \lb{cmf01}

\bibitem{cmf04}
Chen, M.F. (2004).
{\it From Markov Chains to
Non-equilibrium Particle Systems}.
World Scientific. 2$^{\text{nd}}$ ed. (1$^{\text{st}}$ ed., 1992),
Singapore. \lb{cmf04}

\bibitem{cmf05}
Chen, M.F. (2005).
{\it Eigenvalues, Inequalities, and Ergodic Theory}.
Springer, London.  \label{cmf05}

\bibitem{cmf10}
Chen, M.F. (2010).
Speed of stability for birth--death processes.
Front. Math. China 5(3), 379--515.\lb{cmf10}

\bibitem{cmf14}
Chen, M.F. (2014).
{\it Criteria for discrete spectrum of 1D operators}.
\newblock Commu. Math. Stat. 2, 279--309.\label{cmf14}

\bibitem{cmf16}
Chen, M.F. (2016).
{\it Unified speed estimation of various stabilities}.
\newblock Chin. J. Appl. Probab. Statis. 32(1), 1--22.\label{cmf16}

\bibitem{chzhx}
Chen, M.F. and Zhang, X. (2014)
{\it Isospectral operators}.
Commu Math Stat 2, 17--32. \lb{chzhx}

\bibitem{chzh}
Chen, M.F. and Zhang, Y.H. (2014).
{\it Unified representation of formulas for single birth processes}.
Front. Math. China 9(4), 761--796. \lb{chzh}

\bibitem{gl13}
Golub, G.H. and van Loan, C.F. (2013).
{\it Matrix Computations}, 4ed.
Johns Hopkins University Press, Baltimore.\lb{gl13}

\bibitem{hua84}
Hua, L.K. (1984).
{\it Mathematical theory of global optimization on planned economy, \text{\rm (II) and (III)\,(In Chinese)}.}
Kexue Tongbao 13, 769-772. \lb{hua84}

\bibitem{lan1}
Langville, A.N., Meyer, C. D. (2006).
{\it Google's PageRank and Beyond: The Science of Search Engine Rankings.}
Princeton University Press.\lb{lan1}

\bibitem{mc00}
Meyer, C. (2000).
{\it Matrix Analysis and Applied Linear Algebra}.
SIAM.\lb{mc00}

\bibitem{mpg29}
R. von Mises and H. Pollaczek-Geiringer (1929).
{\it Praktische Verfahren der Gleichungsauf\"osung}.
ZAMM - Zeitschrift f\"ur Angewandte Mathematik und Mechanik 9, 152--164.
\lb{mpg29}

\bibitem{hw44}
H. Wielandt (1944).
{\it Beitr\"age zur mathematischen Behandlung komplexer Eigenwertprobleme}.
Teil V: Bestimmung h\"oherer Eigenwerte durch gebrochene Iteration. Bericht B 44/J/37, Aerodynamische Versuchsanstalt G\"ottingen, Germany.\lb{hw44}

\bibitem{jhw65}
Wilkinson, J.H. (1965).
{\it The Algebraic Eigenvalue Problem}.
Oxford University Press, Oxford.  \lb{jhw65}
}
\end{thebibliography}
\end{document}